\documentclass[11pt]{article}
\usepackage[margin=1in]{geometry}
\usepackage{amsmath,amssymb,amsthm,mathtools}
\usepackage{enumitem}
\usepackage{hyperref}
\usepackage{tikz-cd}
\usepackage{graphicx}

\theoremstyle{plain}
\newtheorem{theorem}{Theorem}[section]
\newtheorem{lemma}[theorem]{Lemma}
\newtheorem{proposition}[theorem]{Proposition}
\newtheorem{corollary}[theorem]{Corollary}
\newtheorem*{theoremA}{Theorem A}
\newtheorem*{theoremB}{Theorem B}

\theoremstyle{definition}
\newtheorem{definition}[theorem]{Definition}
\newtheorem{example}[theorem]{Example}
\newtheorem{remark}[theorem]{Remark}

\newcommand{\kk}{k}
\newcommand{\ZZ}{\mathbb{Z}}
\newcommand{\NN}{\mathbb{N}}
\newcommand{\RR}{\mathbb{R}}
\newcommand{\pp}{\mathfrak{p}}
\newcommand{\mm}{\mathfrak{m}}
\newcommand{\link}{\operatorname{lk}}
\newcommand{\cost}{\operatorname{cost}}
\newcommand{\Ext}{\operatorname{Ext}}
\newcommand{\Tor}{\operatorname{Tor}}
\newcommand{\rank}{\operatorname{rank}}
\newcommand{\rH}{\widetilde{H}}
\newcommand{\rC}{\widetilde{C}}
\newcommand{\Hloc}[2]{H^{#1}_{#2}}
\newcommand{\one}{\mathbf 1}
\newcommand{\Nf}{\mathcal N}
\newcommand{\Lc}{\mathcal L}
\newcommand{\Kc}{\mathcal K}

\newcommand{\supp}{\operatorname{supp}}
\newcommand{\res}{\operatorname{res}}

\title{Persistent Bass Numbers and Persistent Lyubeznik Numbers of Stanley--Reisner Rings}

\author{Kaiyue He\thanks{University of Georgia \texttt{hekaiyue@hotmail.com}; corresponding author}
  \and Faisal Suwayyid\thanks{King Fahd University of Petroleum and Minerals \texttt{faisal.suwayyid@kfupm.edu.sa}}
  \and Guo-Wei Wei\thanks{University of Georgia \texttt{Guowei.Wei@uga.edu}}}
\date{}

\begin{document}
\maketitle

\begin{abstract}
We define persistent Lyubeznik numbers and persistent Bass numbers for Stanley–Reisner rings along a simplicial filtration. For Lyubeznik numbers, we show that the  \`Alvarez Montaner–Vahidi complex is natural with respect to inclusions of simplicial complexes, so that persistent Lyubeznik numbers are given by ranks of induced maps on cohomology. We introduce an equivalent costar complex whose terms are reduced homology groups of links of faces and whose differentials are induced by quotient maps. This construction is computationally efficient, requiring neither Alexander duals nor free resolutions. For Bass numbers at coordinate primes, we prove that the Hochster-type formula is natural with respect to inclusions, identifying persistent Bass numbers with the persistent graded Betti numbers of Suwayyid and Wei for filtered links. Both families admit finite barcodes and satisfy bottleneck stability. Computational experiments confirm that the formula-based approach is consistently faster than free resolution methods in Macaulay2, often by one to two orders of magnitude. Examples, including two molecular isomers with identical persistent homology, demonstrate that Lyubeznik barcodes can capture topological information that traditional persistent homology does not detect.


\end{abstract}

\tableofcontents

\section{Introduction}\label{sec:intro}

Persistent Stanley--Reisner theory \cite{SuwayyidWeiPSRT,SuwayyidWeiGraphs} applies invariants of
the face ring $\kk[\Delta]$, such as graded Betti numbers and $f$- and $h$-vectors, to each complex
$\Delta^t$ of a filtration and records how they change. These invariants come from free resolutions.
In \cite{HeSuwayyidWei}, we studied local cohomology at vertex primes along a filtration, which gives
invariants attached to single vertices. The last section of that paper named Bass numbers and
Lyubeznik numbers as the next invariants for further study. This paper treats them.

Throughout, $\Delta$ is a simplicial complex on $[n]=\{1,\dots,n\}$, $f\colon\Delta\to\RR$ is
monotone, and $\Delta^t=\{\sigma\in\Delta:f(\sigma)\le t\}$. For $s\le t$ the inclusion
$\Delta^s\subseteq\Delta^t$ gives a surjection $\kk[\Delta^t]\to\kk[\Delta^s]$ and an inclusion of
ideals $I_{\Delta^t}\subseteq I_{\Delta^s}$. Bass numbers are functorial in the module, so the first
map gives maps from level $t$ to level $s$. Lyubeznik numbers are Bass numbers of local cohomology
modules $\Hloc{i}{I}(S)$, and $I_{\Delta^t}\subseteq I_{\Delta^s}$ gives maps
$\Hloc{i}{I_{\Delta^s}}(S)\to\Hloc{i}{I_{\Delta^t}}(S)$, from level $s$ to level $t$. So the Bass
persistence modules run against the filtration and the Lyubeznik persistence modules run with it.

\paragraph{Bass numbers.} A set $U\subseteq[n]$ of vertices has two coordinate primes,
$\pp_U=(x_j:j\notin U)$ and $\mathfrak q_U=(x_u:u\in U)$, and $\mathfrak q_U=\pp_{[n]\setminus U}$. The variety of $\pp_U$ is the coordinate subspace spanned by $U$; that of $\mathfrak q_U$ is spanned by the other vertices. For
$U=\{i\}$, $\pp_U$ is the vertex prime of \cite{HeSuwayyidWei}. We state results for $\pp_U$ and translate
them to $\mathfrak q_U$ in Section~\ref{ssec:which-prime}, which also explains which prime to use. We use Bass numbers of $\kk[\Delta]$ as an $S$-module; Remark~\ref{rem:over-S} explains this choice.

\begin{theoremA}[Theorems~\ref{thm:bass-static} and \ref{thm:bass-persistent}, Corollaries~\ref{cor:bass-betti-link}, \ref{cor:bass-betti} and \ref{cor:q-prime}, Proposition~\ref{prop:bass-bars}]
The Bass numbers at $\pp_U$ vanish unless $U\in\Delta$, and then $\mu^p(\pp_U,\kk[\Delta])=\beta_{\,n-|U|-p}(\kk[\link_\Delta U])$, which splits as a sum over $G\subseteq[n]\setminus U$
of
\[
\mu^p_{U,G}(\kk[\Delta])=\dim_\kk\rH^{\,p-|G|-1}\bigl((\link_\Delta U)_{[n]\setminus(U\cup G)};\kk\bigr).
\]
This identification is natural for inclusions of complexes. Hence the persistent Bass number
$\mu^{p,s\to t}_{U,G}$ is the rank of restriction
$\rH^{p-|G|-1}((\link_{\Delta^t}U)_W)\to\rH^{p-|G|-1}((\link_{\Delta^s}U)_W)$, $W=[n]\setminus(U\cup G)$.
For $U=\varnothing$ these are the persistent graded Betti numbers of \cite{SuwayyidWeiGraphs}. In degree
$p=0$, the barcode at $\pp_U$ is the single interval on which $U$ is a facet of $\Delta^t$. The Bass numbers
at $\mathfrak q_U$ are those at $\pp_{[n]\setminus U}$: they vanish unless the complement of $U$ is a face, and then they are Betti numbers of $\kk[\link_\Delta([n]\setminus U)]$, a complex on $U$.
\end{theoremA}

For $U=\varnothing$, Theorem~A is equivalent, through the Koszul self-duality of
Remark~\ref{rem:koszul-duality} (homological degree $n-p$, multidegree
$\one_{[n]\setminus G}$), to the persistent Hochster formula
\cite[Thm~3.2]{SuwayyidWeiGraphs}. For general $U$ it is the same formula for the
filtered link $(\link_{\Delta^t}U)_t$ over $S_{U^c}$, obtained by the localization of
Lemma~\ref{lem:star-loc} and Proposition~\ref{prop:bass-reduction}
(Corollary~\ref{cor:bass-betti}). So Theorem~A gives no new invariant of the link. What
we add is a complete proof, with signs, that the Hochster isomorphism commutes with the
maps induced by inclusions (Theorem~\ref{thm:bass-static}(b)); for this step
\cite{SuwayyidWeiGraphs} gives an idea of proof.

\paragraph{Lyubeznik numbers.} Let $\Delta^\vee=\{H\subseteq[n]:[n]\setminus H\notin\Delta\}$ be the Alexander
dual. For $J\in\Delta$ the restriction $(\Delta^\vee)_{[n]\setminus J}$ is the Alexander dual of $\link_\Delta J$
inside the ground set $[n]\setminus J$ (Lemma~\ref{lem:dual-of-link}).

\begin{theoremB}[Theorems~\ref{thm:lyub-complex}, \ref{thm:lyub-persistent},
\ref{thm:lyub-stability} and \ref{thm:link-model}, Corollary~\ref{cor:costar-persistent}]
Let $\Delta\neq2^{[n]}$. There is an isomorphism
$\Ext^p_S(\kk,\Hloc{n-q}{I_\Delta}(S))\cong H^p(\Lc^\bullet_q(\Delta))$, where
$\Lc^p_q(\Delta)=\bigoplus_{|J|=p}\rH^{\,n-q-2}((\Delta^\vee)_{[n]\setminus J};\kk)$ with
signed restriction maps, and it is natural for inclusions $\Delta'\subseteq\Delta''$ of
complexes different from $2^{[n]}$. Hence, for levels $\Delta^s\subseteq\Delta^t$
different from $2^{[n]}$,
\[
\lambda^{s\to t}_{p,q}=\rank\Bigl(H^p(\Lc^\bullet_q(\Delta^s))\to H^p(\Lc^\bullet_q(\Delta^t))\Bigr).
\]
The complex $\Lc^\bullet_q(\Delta)$ is naturally isomorphic to the costar complex
$\Kc^\bullet_q(\Delta)$, whose $p$-th term is
$\bigoplus_{F\in\Delta,\,|F|=p}\rH_{q-p-1}(\link_\Delta F;\kk)$; the same rank formula holds
with $\Kc^\bullet_q$, also when $\Delta^t=2^{[n]}$. Each Lyubeznik persistence module has a
finite barcode $\mathcal B_{p,q}(f)$, and
$d_b(\mathcal B_{p,q}(f),\mathcal B_{p,q}(g))\le\|f-g\|_\infty$.
\end{theoremB}

The equality $\lambda_{p,q}(\kk[\Delta])=\dim_\kk H^p(\Lc^\bullet_q(\Delta))$ for a single
complex is due to \`Alvarez Montaner and Vahidi \cite{AlvarezMontanerVahidi}, building on
Musta\c{t}\u{a} \cite{Mustata}. What is new in Theorem~B is the naturality in $\Delta$,
the persistent invariants it gives, and the costar model.

\begin{remark}

\textbf{What is known.} The static formulas are close to known results. Bass numbers of $\kk[\Delta]$
at coordinate primes follow from Hochster's formula \cite{Hochster} by localization, and related statements
appear in the theory of squarefree modules \cite{Yanagawa2000,MillerSturmfels}. The graded pieces of
$\Hloc{i}{I_\Delta}(S)$ were computed by Mustață \cite{Mustata}. Gr\"abe \cite{Grabe1984} described the multigraded structure of local cohomology of Stanley--Reisner rings in terms of relative simplicial cohomology of costars, including the corresponding multiplication maps by the variables. He also established the link--costar homology isomorphism used below. Combinatorial descriptions of Lyubeznik numbers of Stanley--Reisner rings were given by Àlvarez Montaner and Vahidi \cite{AlvarezMontanerVahidi} through linear strands of the Alexander dual ideal, and by Àlvarez Montaner and Yanagawa \cite{AlvarezMontanerYanagawa}, who also prove that $\lambda_{p,q}(\kk[\Delta])$ depends only on $|\Delta|$ and $\operatorname{char}\kk$ \cite[Theorem~5.3]{AlvarezMontanerYanagawa}. In their proof $\lambda_{p,q}$ is computed from the sheaf of local homology of $|\Delta|$, whose stalks are the groups $\rH_{q-|J|-1}(\link_\Delta J)$ that appear in Lemma~\ref{lem:costar-link}. Yanagawa computed Bass numbers of these local cohomology modules
\cite{Yanagawa2001}.

\medskip\noindent
\textbf{What is new.} What we add is the following: (i) versions of these descriptions in which the identifications are natural for inclusions of complexes, including the multiplication maps by the variables, which is what persistence needs; this includes the costar and link complex of Theorem~\ref{thm:link-model}, a chain-level form of the description in \cite{AlvarezMontanerYanagawa}, building on the classical link--costar identification and multiplication maps of \cite{Grabe1984}, with explicit signs; (ii) the persistent invariants and their finiteness and stability;  (iii) statements about bars, such as facet lifetimes for $\mu^0$ and the left endpoints of $\lambda_{q,q}$ bars; and (iv) a worked example.

\medskip\noindent

\end{remark}

\paragraph{Paper Outline.} Section~\ref{sec:background} fixes conventions. Sections~\ref{sec:bass-static}
and~\ref{sec:bass-persistent} give basics of the Bass numbers; Section~\ref{ssec:which-prime} discusses the choice between
$\pp_U$ and $\mathfrak q_U$. Section~\ref{sec:lyub-static} gives the basics of the Lyubeznik complex,
Section~\ref{sec:lyub-persistent} its persistent version, and Section~\ref{sec:link-model} gives an equivalent
costar model and its interpretation through links. Example~\ref{ex:additional-lyubeznik-bars} is a first example, and
Appendix~\ref{app:computations} computes all the invariants for two further examples and
compares the  two computations.

\section{Conventions and background}\label{sec:background}

\subsection{Simplicial complexes}

A simplicial complex $K$ on a finite ground set $V$ is a family of subsets of $V$ closed under taking
subsets. Elements of $V$ need not be vertices of $K$. The \emph{void complex} has no faces, and
$\{\varnothing\}$ is the complex whose only face is $\varnothing$. Reduced cochains $\rC^\bullet(K;\kk)$
have basis the faces, a face with $m$ elements sitting in degree $m-1$; with the elements of $V$ ordered,
the coboundary is $\delta(\tau^\vee)=\sum_{j}(-1)^{\nu_\tau(j)}(\tau\cup\{j\})^\vee$, the sum over
$j\notin\tau$ with $\tau\cup\{j\}\in K$, where $\nu_A(j)=|\{u\in A:u<j\}|$. Thus
$\rH^{-1}(\{\varnothing\};\kk)=\kk$, the void complex has zero reduced cohomology, and a simplex $2^V$
with $V\neq\varnothing$ has zero reduced cohomology. Over a field, $\rH^d$ and $\rH_d$ have equal
dimension. For $W\subseteq V$, $K_W=\{\sigma\in K:\sigma\subseteq W\}$. For $F\subseteq V$,
\[
\link_K F=\{\tau\subseteq V\setminus F:\tau\cup F\in K\},
\]
a complex on $V\setminus F$; it is void if $F\notin K$. The \emph{costar} of $F$ is
\[
\cost_K F=\{\tau\in K:F\not\subseteq\tau\}.
\]
It is a subcomplex of $K$. If $F\notin K$, then
$\cost_K F=K$, while $\cost_K\varnothing$ is void.
For $K\neq 2^V$ the \emph{Alexander dual} of $K$ in
$V$ is $K^{*V}=\{H\subseteq V:V\setminus H\notin K\}$. We write $\Delta^\vee=\Delta^{*[n]}$.

\begin{theorem}[Combinatorial Alexander duality {\cite{BjornerTancer}}]\label{thm:alexander}
Let $V\neq\varnothing$ and let $K\neq 2^V$ be a simplicial complex on $V$. Then
$\rH^i(K;\kk)\cong\rH_{|V|-i-3}(K^{*V};\kk)$ for all $i$.
\end{theorem}

\begin{lemma}\label{lem:dual-of-link}
Let $\Delta\neq2^{[n]}$ and $F\subsetneq[n]$. If $F\in\Delta$ then $(\Delta^\vee)_{[n]\setminus F}$ is the
Alexander dual of $\link_\Delta F$ in $[n]\setminus F$. If $F\notin\Delta$ then
$(\Delta^\vee)_{[n]\setminus F}=2^{[n]\setminus F}$. In both cases $(\Delta^\vee)_{[n]\setminus F}$ contains
$\varnothing$ and depends only on $\link_\Delta F$.
\end{lemma}

\begin{proof}
Put $V=[n]\setminus F$. For $H\subseteq V$ we have $[n]\setminus H=F\cup(V\setminus H)$, so $H\in\Delta^\vee$
if and only if $V\setminus H\notin\link_\Delta F$. If $F\in\Delta$ this is the definition of the Alexander
dual of $\link_\Delta F$ in $V$. If $F\notin\Delta$ the link is void and every $H\subseteq V$ qualifies.
Since $[n]\notin\Delta$, $\varnothing\in\Delta^\vee$.
\end{proof}

\subsection{Koszul cochain complexes}

Let $V$ be a finite ordered set, $S_V=\kk[x_v:v\in V]$ with the $\ZZ^V$-grading, $\mm_V$ its graded
maximal ideal, and $N$ a $\ZZ^V$-graded $S_V$-module. For $J\subseteq V$ write $\one_J$ for its indicator
vector and $x_J=\prod_{j\in J}x_j$. The Koszul cochain complex $K^\bullet(N)$ has
$K^p(N)=\bigoplus_{|J|=p}N\,e_J$, where $u\,e_J$ has degree $\deg u-\one_J$, and
\[
d(u\,e_J)=\sum_{j\in V\setminus J}(-1)^{\nu_J(j)}x_ju\,e_{J\cup\{j\}}.
\]
It is $\operatorname{Hom}_{S_V}(K_\bullet,N)$ for the graded Koszul resolution $K_\bullet$ of $\kk$, so
$H^p(K^\bullet(N))\cong\Ext^p_{S_V}(\kk,N)$ as graded vector spaces, naturally in $N$ \cite[\S1.6]{BrunsHerzog}.
Its degree-$a$ part is $K^p(N)_a=\bigoplus_{|J|=p}N_{a+\one_J}e_J$.

\begin{lemma}\label{lem:cone}
Let $a\in\ZZ^V$ and $j\in V$. If multiplication $x_j\colon N_{a+\one_J}\to N_{a+\one_J+\mathbf e_j}$ is bijective
for every $J\subseteq V\setminus\{j\}$, then $K^\bullet(N)_a$ is acyclic.
\end{lemma}

\begin{proof}
Let $A\subseteq K^\bullet(N)_a$ be spanned by the $u\,e_J$ with $j\in J$ and $Q$ the span of those with
$j\notin J$. The differential only enlarges $J$, so $A$ is a subcomplex with quotient $Q$. Write the
differential as $d_Q$ on $Q$, $d_A$ on $A$, and $\varphi\colon Q\to A$, $u\,e_J\mapsto(-1)^{\nu_J(j)}x_ju\,e_{J\cup\{j\}}$.
From $d^2=0$ we get $d_A\varphi=-\varphi d_Q$, and $\varphi$ is bijective in each degree by hypothesis.
The connecting map of $0\to A\to K^\bullet(N)_a\to Q\to0$ is induced by $\varphi$: a cocycle $q$ of $Q$ lifts to
$q$, and $d q=\varphi(q)$. Since $\varphi$ is a bijection that anticommutes with the differentials, it
induces isomorphisms $H^p(Q)\to H^{p+1}(A)$. The long exact sequence then gives $H(K^\bullet(N)_a)=0$.
\end{proof}

\subsection{Bass numbers and Lyubeznik numbers}

For a prime $\pp$ of $S$ and an $S$-module $M$, the Bass numbers are
$\mu^p(\pp,M)=\dim_{\kappa(\pp)}\Ext^p_{S_\pp}(\kappa(\pp),M_\pp)$, where $\kappa(\pp)=S_\pp/\pp S_\pp$.
They count copies of the injective hull of $S/\pp$ in a minimal injective resolution of $M$
\cite[\S3.2]{BrunsHerzog}.

\begin{remark}\label{rem:over-S}
We always take Bass numbers of $\kk[\Delta]$ as an $S$-module. Since $S$ is regular, they vanish for
$p>n$. Bass numbers of $\kk[\Delta]$ as a module over itself are in general nonzero in infinitely many
degrees, because a local ring has finite injective dimension over itself only when it is Gorenstein.
Lyubeznik numbers are also defined over the regular ring, so working over $S$ treats both invariants
in the same way.
\end{remark}

Let $I\subseteq S$ be a homogeneous ideal and $R=S/I$. The Lyubeznik numbers of $R$ are
\[
\lambda_{p,q}(R)=\mu^p\bigl(\mm,\Hloc{n-q}{I}(S)\bigr),
\]
which are finite and depend only on $R$ \cite{Lyubeznik,HunekeSharp}. Here Lyubeznik numbers of the
graded ring $R$ mean those of its localization at $\mm$. Local cohomology commutes with localization,
and $\Ext^p_S(\kk,N)$ is annihilated by $\mm$, so
$\lambda_{p,q}(R)=\dim_\kk\Ext^p_S(\kk,\Hloc{n-q}{I}(S))$.

\subsection{Filtrations and persistence modules}

A function $f\colon\Delta\to\RR$ is monotone if $f(\tau)\le f(\sigma)$ for $\tau\subseteq\sigma$. It gives
the filtration $\Delta^t=\{\sigma\in\Delta:f(\sigma)\le t\}$ by subcomplexes on $[n]$. For
$t<f(\varnothing)$ the complex $\Delta^t$ is void and $\kk[\Delta^t]=0$. We keep the standing convention of
\cite{HeSuwayyidWei} that every singleton is a face of $\Delta$; the complexes $\Delta^t$ may miss vertices.

An $\RR$ or $\RR^{op}$-indexed persistence module is a functor from $(\RR,\le)$ to $\kk$-vector spaces; an
$\RR^{\mathrm{op}}$-indexed module has its maps going from $t$ to $s$ for $s\le t$ \cite{BubenikScott}. A
module is of \textit{finite type} if it is pointwise finite-dimensional and there is a finite set $C\subseteq\RR$
such that all maps between values in the same component of $\RR\setminus C$, or from a point of $C$ into
the component to its right, are isomorphisms. Such a module is a finite direct sum of interval modules,
uniquely up to order \cite{CrawleyBoevey}; the multiset of intervals is its barcode. Substituting
$t\mapsto-t$ turns an $\RR^{\mathrm{op}}$-indexed module into an $\RR$-indexed one and reflects its barcode,
which does not change bottleneck distances. For $\RR$-indexed modules we use $\delta$-interleavings and the
isometry theorem of Chazal, de Silva, Glisse and Oudot \cite{ChazalDeSilvaGlisseOudot}: $\delta$-interleaved
modules of finite type have barcodes at bottleneck distance at most $\delta$.

\begin{lemma}\label{lem:functor-stability}
Let $\Phi$ be a functor from the poset of subcomplexes of $\Delta$ (ordered by inclusion) to finite-dimensional
vector spaces, covariant or contravariant. For monotone $f$ let $\Phi(f)$ be the persistence module
$t\mapsto\Phi(\Delta^t_f)$. Then $\Phi(f)$ is of finite type with critical values in $\{f(\sigma):\sigma\in\Delta\}$,
and the barcodes of $\Phi(f)$ and $\Phi(g)$ are at bottleneck distance at most $\|f-g\|_\infty$.
\end{lemma}

\begin{proof}
If no value of $f$ lies in $(s,t]$ then $\Delta^s_f=\Delta^t_f$ and the map is the identity; this gives
finite type. Put $\delta=\|f-g\|_\infty$. If $f(\sigma)\le t$ then $g(\sigma)\le t+\delta$, so
$\Delta^t_f\subseteq\Delta^{t+\delta}_g\subseteq\Delta^{t+2\delta}_f$, and symmetrically. Applying $\Phi$ to
these inclusions gives maps in both directions whose composites are the structure maps over $2\delta$,
because $\Phi$ is a functor on a poset. In the contravariant case, first substitute $t\mapsto-t$. This is
a $\delta$-interleaving, and the isometry theorem gives the bound.
\end{proof}

\section{Bass numbers at the primes of a vertex set}\label{sec:bass-static}

Throughout this section $U\subseteq[n]$ is a set of vertices and $U^c=[n]\setminus U$. Two primes of $S$ are
attached to $U$:
\[
\pp_U=(x_j:j\in U^c),\qquad \mathfrak q_U=(x_u:u\in U).
\]
Since $S/\pp_U\cong\kk[x_u:u\in U]$, the variety of $\pp_U$ is the coordinate subspace $\kk^U$ spanned by the
chosen vertices, and the variety of $\mathfrak q_U$ is the complementary coordinate subspace $\kk^{U^c}$. The two
families are the same family indexed differently: $\mathfrak q_U=\pp_{U^c}$. For $U=\{i\}$, $\pp_U$ is the vertex
prime of \cite{HeSuwayyidWei} and $\mathfrak q_U=(x_i)$ is the prime of the coordinate hyperplane. We develop the
theory for $\pp_U$ with $U$ arbitrary; Section~\ref{ssec:which-prime} translates it to $\mathfrak q_U$ and explains
which prime to use. Multidegrees in $\ZZ^n$ are written $(c;b)$ with $c\in\ZZ^U$ and $b\in\ZZ^{U^c}$.

\begin{lemma}\label{lem:support}
$\pp_U\supseteq I_\Delta$ if and only if $U\in\Delta$, and then $\kk[\Delta]/\pp_U\kk[\Delta]\cong\kk[x_u:u\in U]$.
Hence $\pp_U$ lies in the support of $\kk[\Delta]$ exactly when $U$ is a face, and $\kk[\Delta]_{\pp_U}=0$ when
$U\notin\Delta$. In particular $\mathfrak q_U$ lies in the support of $\kk[\Delta]$ exactly when $U^c\in\Delta$.
\end{lemma}

\begin{proof}
A monomial $x_N$ lies in $\pp_U$ if and only if $N\not\subseteq U$. So $I_\Delta\subseteq\pp_U$ if and only if no
nonface is contained in $U$, which says $U\in\Delta$. If $U\notin\Delta$ then $x_U\in I_\Delta$ while
$x_U\notin\pp_U$, so $x_U$ is a unit in $\kk[\Delta]_{\pp_U}$ and also zero. The support of a finitely generated
module consists of the primes containing its annihilator.
\end{proof}

\subsection{Reduction to the maximal ideal of a smaller ring}

\begin{definition}
For a $\ZZ^n$-graded $S$-module $M$ let $M_{x_U}$ be its localization at the powers of $x_U=\prod_{u\in U}x_u$, and let
\[
M_{(U)}=\bigoplus_{b\in\ZZ^{U^c}}(M_{x_U})_{(0;b)},
\]
the part of $M_{x_U}$ of degree $0$ in the $U$-coordinates. It is a $\ZZ^{U^c}$-graded $S_{U^c}$-module.
\end{definition}

\begin{lemma}\label{lem:star-loc}
\begin{enumerate}[label=\textup{(\alph*)},nosep]
\item $M_{x_U}\cong M_{(U)}\otimes_\kk\kk[x_u^{\pm1}:u\in U]$ as $\ZZ^n$-graded modules over
$S_{x_U}=S_{U^c}\otimes_\kk\kk[x_u^{\pm1}:u\in U]$, naturally in $M$.
\item If $U\in\Delta$, the monomials of $\kk[\Delta]_{(U)}$ are the $x^b$ with $b\in\NN^{U^c}$ and
$\supp b\in\link_\Delta U$, so $\kk[\Delta]_{(U)}=\kk[\link_\Delta U]$ as $\ZZ^{U^c}$-graded $S_{U^c}$-modules.
If $U\notin\Delta$, $\kk[\Delta]_{(U)}=0$.
\item For $\Delta'\subseteq\Delta$ the map $\kk[\Delta]_{(U)}\to\kk[\Delta']_{(U)}$ induced by
$\kk[\Delta]\to\kk[\Delta']$ is the quotient $\kk[\link_\Delta U]\to\kk[\link_{\Delta'}U]$, which sends each
monomial to itself or to $0$.
\end{enumerate}
\end{lemma}

\begin{proof}
(a) For $c\in\ZZ^U$, multiplication by the unit $x^c$ is a bijection $(M_{x_U})_{(0;b)}\to(M_{x_U})_{(c;b)}$
that commutes with multiplication by $S_{U^c}$. Summing over $c$ gives the isomorphism.
(b) By the localization lemma of \cite{HeSuwayyidWei} (see also \cite[Lemma~5.3.6]{BrunsHerzog}), $x^a$ is nonzero in
$\kk[\Delta]_{x_U}$ exactly when $a_j\ge0$ for $j\notin U$ and $\supp_+(a)\cup U\in\Delta$; each such
graded piece is spanned by $x^a$. For $a=(0;b)$ the condition says $b\ge0$ and $\supp b\cup U\in\Delta$, and
multiplication by $x_j$, $j\notin U$, obeys the same rule on both sides. If $U\notin\Delta$ no monomial qualifies.
(c) The surjection sends each monomial to itself or to $0$, and so does its localization.
\end{proof}

\begin{proposition}\label{prop:bass-reduction}
Let $M$ be a $\ZZ^n$-graded $S$-module and $U\subseteq[n]$. For every $p$,
\[
\mu^p(\pp_U,M)=\dim_\kk\Ext^p_{S_{U^c}}\bigl(\kk,M_{(U)}\bigr).
\]
If $M\to M'$ is a degree-preserving map with finite Bass numbers on both sides, the induced maps
$\Ext^p_{S_{\pp_U}}(\kappa(\pp_U),M_{\pp_U})\to\Ext^p_{S_{\pp_U}}(\kappa(\pp_U),M'_{\pp_U})$ and
$\Ext^p_{S_{U^c}}(\kk,M_{(U)})\to\Ext^p_{S_{U^c}}(\kk,M'_{(U)})$ have the same rank.
\end{proposition}

\begin{proof}
Since $x_u\notin\pp_U$ for $u\in U$, $M_{\pp_U}$ is a localization of $M_{x_U}$. The ring map
$S_{U^c}\to S_{\pp_U}$ is flat, being the free extension $S_{U^c}\to S_{x_U}$ followed by a localization.
By Lemma~\ref{lem:star-loc}(a), $M_{\pp_U}\cong M_{(U)}\otimes_{S_{U^c}}S_{\pp_U}$, and
$\kk\otimes_{S_{U^c}}S_{\pp_U}=S_{\pp_U}/\pp_US_{\pp_U}=\kappa(\pp_U)$. Since $\kk$ has a resolution by finite
free $S_{U^c}$-modules, flat base change gives
\[
\Ext^p_{S_{U^c}}(\kk,M_{(U)})\otimes_{S_{U^c}}S_{\pp_U}\cong\Ext^p_{S_{\pp_U}}(\kappa(\pp_U),M_{\pp_U}),
\]
naturally in $M$. The left side is $E\otimes_\kk\kappa(\pp_U)$ with $E=\Ext^p_{S_{U^c}}(\kk,M_{(U)})$, because
$E$ is annihilated by $\mm_{U^c}$. Extending scalars from $\kk$ to $\kappa(\pp_U)$ preserves dimensions and
ranks of linear maps.
\end{proof}

\begin{definition}\label{def:graded-bass}
For $b\in\ZZ^{U^c}$ put $\mu^p_U(M)_b=\dim_\kk\Ext^p_{S_{U^c}}(\kk,M_{(U)})_b$, and for $G\subseteq U^c$ put
$\mu^p_{U,G}(M)=\mu^p_U(M)_{-\one_G}$.
\end{definition}

\begin{remark}\label{rem:koszul-duality}
If $I^\bullet$ is a minimal $\ZZ^{U^c}$-graded injective resolution of $M_{(U)}$, the complex
$\operatorname{Hom}(\kk,I^\bullet)$ has zero differential, so $\mu^p_U(M)_b$ is the dimension of the degree-$b$
part of the socle of $I^p$. For $\ZZ^n$-graded injective resolutions see \cite{GotoWatanabe} and
\cite[Ch.~11]{MillerSturmfels}. Over a polynomial ring $S_V$ in $m$ variables the Koszul complex resolving $\kk$
is self-dual, which gives $\Ext^p_{S_V}(\kk,N)_b\cong\Tor^{S_V}_{m-p}(\kk,N)_{b+\one_V}$ for every graded $N$.
So Bass numbers at the maximal ideal are Betti numbers: $\mu^p(\mm_V,N)=\beta_{m-p}(N)$.
\end{remark}

\subsection{A Hochster-type formula}

\begin{theorem}\label{thm:bass-static}
Let $\Gamma$ be a simplicial complex on a finite ordered set $V$ and $p\ge0$.
\begin{enumerate}[label=\textup{(\alph*)},nosep]
\item $\Ext^p_{S_V}(\kk,\kk[\Gamma])_b=0$ unless $b=-\one_G$ for some $G\subseteq V$.
\item For $G\subseteq V$ there is an isomorphism
\[
\Ext^p_{S_V}(\kk,\kk[\Gamma])_{-\one_G}\cong\rH^{\,p-|G|-1}\bigl(\Gamma_{V\setminus G};\kk\bigr).
\]
For a subcomplex $\Gamma'\subseteq\Gamma$, the map induced by $\kk[\Gamma]\to\kk[\Gamma']$ corresponds to
restriction $\rH^{p-|G|-1}(\Gamma_{V\setminus G})\to\rH^{p-|G|-1}(\Gamma'_{V\setminus G})$.
\end{enumerate}
\end{theorem}

\begin{proof}
We use $K^\bullet(\kk[\Gamma])_b=\bigoplus_J\kk[\Gamma]_{b+\one_J}e_J$. Recall that $\kk[\Gamma]_a$ is
$\kk\cdot x^a$ if $a\ge0$ and $\supp a\in\Gamma$, and $0$ otherwise.

(a) If $b_j\le-2$ for some $j$, every term is zero. If $b_j\ge1$ for some $j$, then for $J\not\ni j$ the
$j$-th coordinate of $b+\one_J$ is at least $1$, so multiplication by $x_j$ does not change the support and
is bijective from degree $b+\one_J$ to $b+\one_J+\mathbf e_j$. Lemma~\ref{lem:cone} gives acyclicity. The
remaining case is $b\in\{0,-1\}^V$.

(b) Let $b=-\one_G$. The summand for $J$ is nonzero exactly when $G\subseteq J$ and $\tau:=J\setminus G\in\Gamma$.
Then $\tau\subseteq V\setminus G$ and $|\tau|=p-|G|$. Let $\varepsilon(\tau)=(-1)^{|\{(g,u)\in G\times\tau:g<u\}|}$
and define $\Psi(x^{\one_\tau}e_{G\cup\tau})=\varepsilon(\tau)\,\tau^\vee\in\rC^{\,p-|G|-1}(\Gamma_{V\setminus G})$.
This is a bijection of bases. For $j\notin G\cup\tau$ we have $\nu_{G\cup\tau}(j)=\nu_G(j)+\nu_\tau(j)$ and
$\varepsilon(\tau\cup\{j\})=\varepsilon(\tau)(-1)^{\nu_G(j)}$, and $x_jx^{\one_\tau}$ is $x^{\one_{\tau\cup\{j\}}}$
if $\tau\cup\{j\}\in\Gamma$ and $0$ otherwise. Hence
$\Psi d(x^{\one_\tau}e_{G\cup\tau})=\varepsilon(\tau)\sum_j(-1)^{\nu_\tau(j)}(\tau\cup\{j\})^\vee=\delta\Psi(x^{\one_\tau}e_{G\cup\tau})$,
so $\Psi$ is an isomorphism of cochain complexes. The same sign computation appears in the proof of Hochster's
formula in \cite{HeSuwayyidWei}. The signs do not depend on $\Gamma$, and $\kk[\Gamma]\to\kk[\Gamma']$ sends
$x^{\one_\tau}$ to itself if $\tau\in\Gamma'$ and to $0$ otherwise, which under $\Psi$ is restriction of cochains.
\end{proof}

\begin{corollary}\label{cor:bass-static}
Let $U\subseteq[n]$, $G\subseteq U^c$, and $W=[n]\setminus(U\cup G)$. Then
\[
\mu^p_{U,G}(\kk[\Delta])=
\begin{cases}\dim_\kk\rH^{\,p-|G|-1}\bigl((\link_\Delta U)_W;\kk\bigr),&U\in\Delta,\\ 0,&U\notin\Delta,\end{cases}
\qquad
\mu^p(\pp_U,\kk[\Delta])=\sum_{G\subseteq U^c}\mu^p_{U,G}(\kk[\Delta]).
\]
In particular:
\begin{enumerate}[label=\textup{(\roman*)},nosep]
\item For $U=\varnothing$, $\mu^p_{\varnothing,G}(\kk[\Delta])=\dim\rH^{p-|G|-1}(\Delta_{[n]\setminus G})$, which is
Hochster's formula for the graded Betti number $\beta_{n-p,\one_{[n]\setminus G}}(\kk[\Delta])$.
\item $\mu^0(\pp_U,\kk[\Delta])$ is $1$ if $U$ is a facet of $\Delta$ and $0$ otherwise.
\item $\mu^p(\pp_U,\kk[\Delta])=0$ for $p>n-|U|$, and $\mu^{n-|U|}(\pp_U,\kk[\Delta])=1$ if $U\in\Delta$.
\item For $\Delta=2^{[n]}$, $\mu^p(\pp_U,S)=1$ if $p=n-|U|$ and $0$ otherwise.
\end{enumerate}
\end{corollary}

\begin{proof}
The formula follows from Proposition~\ref{prop:bass-reduction}, Lemma~\ref{lem:star-loc}(b) and
Theorem~\ref{thm:bass-static} with $V=U^c$ and $\Gamma=\link_\Delta U$. For (i) compare with Hochster's
formula $\beta_{i,\one_W}(\kk[\Delta])=\dim\rH^{|W|-i-1}(\Delta_W)$ at $W=[n]\setminus G$, $i=n-p$.
For (ii), $p=0$ forces degree $-|G|-1\ge-1$, so $G=\varnothing$, and $\rH^{-1}(\link_\Delta U)\ne0$ exactly
when $\link_\Delta U=\{\varnothing\}$, that is, when $U$ is a facet. For (iii), the degree $p-|G|-1$ exceeds
$\dim(\link_\Delta U)_W\le|W|-1=n-|U|-|G|-1$ when $p>n-|U|$; for $p=n-|U|$ only $G=U^c$, $W=\varnothing$
contributes, with $\rH^{-1}(\{\varnothing\})=\kk$. For (iv), $\link U$ is the simplex on $U^c$, whose restriction
to $W$ is acyclic unless $W=\varnothing$; then $G=U^c$ and the degree is $-1$ exactly when $p=n-|U|$.
\end{proof}

Item (ii) agrees with the fact that the associated primes of $\kk[\Delta]$ are the $\pp_U$ with $U$ a facet
\cite[\S5.1]{BrunsHerzog}. Item (iii) says that the Bass numbers at $\pp_U$ live in the range
$0\le p\le\operatorname{ht}\pp_U$, and (iv) agrees with $S$ being Gorenstein.

\begin{corollary}\label{cor:bass-betti-link}
For $U\in\Delta$ and $G\subseteq U^c$,
\[
\mu^p(\pp_U,\kk[\Delta])=\beta^{S_{U^c}}_{\,n-|U|-p}\bigl(\kk[\link_\Delta U]\bigr),\qquad
\mu^p_{U,G}(\kk[\Delta])=\beta^{S_{U^c}}_{\,n-|U|-p,\ \one_{[n]\setminus(U\cup G)}}\bigl(\kk[\link_\Delta U]\bigr).
\]
So the Bass numbers of $\kk[\Delta]$ at $\pp_U$ are the graded Betti numbers of the Stanley--Reisner ring of the
link of $U$, over the polynomial ring in the variables outside $U$.
\end{corollary}

\begin{proof}
Proposition~\ref{prop:bass-reduction}, Lemma~\ref{lem:star-loc}(b), and Remark~\ref{rem:koszul-duality} with
$V=U^c$, $m=n-|U|$.
\end{proof}

\begin{remark}[Comparison with vertex-prime local cohomology]\label{rem:vertex-prime}
For $U=\{i\}$, Corollary~\ref{cor:bass-betti-link} shows that the Bass numbers of $\kk[\Delta]$ at $\pp_i$
depend only on $\link_\Delta(i)$. The local cohomology $\Hloc{q}{\pp_i}(\kk[\Delta])$ studied in
\cite{HeSuwayyidWei} splits into a summand for $\link_\Delta(i)$ and a summand for the deletion of $i$, and
the second summand depends on all of $\Delta$ away from $i$. Its graded pieces are reduced cohomology of
links $\link_\Delta(\sigma\cup\{i\})$ and of $\link_\Delta(\sigma)\cap\operatorname{del}_\Delta(i)$, while the
Bass numbers use induced subcomplexes of $\link_\Delta(i)$.
\end{remark}

\subsection{The prime generated by the chosen vertices}\label{ssec:which-prime}

Everything above applies to $\mathfrak q_U=(x_u:u\in U)$ after replacing $U$ by $U^c$.

\begin{corollary}\label{cor:q-prime}
Let $U\subseteq[n]$. Then $\mu^p(\mathfrak q_U,\kk[\Delta])=0$ for all $p$ unless $U^c\in\Delta$, and in that case
\[
\mu^p(\mathfrak q_U,\kk[\Delta])=\beta^{S_U}_{\,|U|-p}\bigl(\kk[\link_\Delta U^c]\bigr)
=\sum_{G\subseteq U}\dim_\kk\rH^{\,p-|G|-1}\bigl((\link_\Delta U^c)_{U\setminus G};\kk\bigr),
\]
where $\link_\Delta U^c$ is a complex on the vertex set $U$. These numbers vanish for $p>|U|$, and
$\mu^0(\mathfrak q_U,\kk[\Delta])=1$ exactly when $U^c$ is a facet.
\end{corollary}

\begin{remark}[Which prime to use]\label{rem:which-prime}
Let $U$ be a set of vertices of interest, for instance a group of atoms in a molecule or a cluster in a point
cloud. The two primes ask different questions.

\emph{$\pp_U$: the geometry near the chosen vertices.} The Bass numbers at $\pp_U$ are the Betti numbers of
$\kk[\link_\Delta U]$, so they record how the rest of $\Delta$ is attached to the face $U$. They are nonzero only
when $U$ is a face (Lemma~\ref{lem:support}), and in a filtration they vanish for $t<f(U)$
(Proposition~\ref{prop:bass-bars}). This is the prime to use when $U$ is small and its members are close,
which is the typical situation for local analysis: $U$ is a single vertex, an edge, or a small simplex, and
the invariant describes its neighborhood.

\emph{$\mathfrak q_U$: the geometry near the complement.} The Bass numbers at $\mathfrak q_U$ are the Betti numbers of
$\kk[\link_\Delta U^c]$, a complex on $U$, so they record how the chosen vertices are attached to the face
$U^c$. They are nonzero only when $U^c$ is a face, that is, when all vertices outside $U$ span a simplex of
$\Delta$. For a complex built from a point cloud with many points and small faces this never happens, and
every Bass number at $\mathfrak q_U$ is zero at every filtration level. The prime $\mathfrak q_U$ is the right choice
when $U^c$ is a face, for example when $\Delta$ is a simplex with a few vertices removed or a cone over a small
complex with the apex set outside $U$, and when the question is how a few marked vertices sit relative to a
large face. In the persistent setting the bars at $\mathfrak q_U$ start no earlier than $f(U^c)$, so they are
visible only after the complement has filled in.

\emph{$U$ not a face.} If $U\notin\Delta$ then $\pp_U$ is outside the support of $\kk[\Delta]$ and there are no
Bass numbers to compute. Two replacements are available: the collection $\{\mu^p(\pp_F,\kk[\Delta]):F\subseteq U,
F\in\Delta\}$ over the faces contained in $U$, which is a face-indexed family in the sense of
\cite{HeSuwayyidWei}; or local cohomology with support in the ideal $(x_j:j\notin U)$, which is defined whether or
not this ideal is prime and is treated in \cite[\S7]{HeSuwayyidWei}. Neither is a Bass number, and the
persistent theory of Section~\ref{sec:bass-persistent} applies to each $\pp_F$ in the first collection separately.

In summary: choose $\pp_U$ to localize at the vertices in $U$, and $\mathfrak q_U$ to localize at the vertices
outside $U$. Writing down $(x_{i_1},\dots,x_{i_r})$ for a set of chosen points, which is the natural first
guess, gives $\mathfrak q_U$ and therefore localizes at everything except the chosen points.
\end{remark}

\section{Persistent Bass numbers}\label{sec:bass-persistent}

\begin{definition}\label{def:pers-bass}
Let $f$ be monotone, $U\subseteq[n]$, $G\subseteq U^c$, $p\ge0$. The \emph{Bass persistence module} is the
$\RR^{\mathrm{op}}$-indexed module
\[
\mathbb M^p_{U,G}(f)=\Bigl(\Ext^p_{S_{U^c}}\bigl(\kk,\kk[\Delta^t]_{(U)}\bigr)_{-\one_G}\Bigr)_{t\in\RR},
\]
with maps from $t$ to $s$ ($s\le t$) induced by $\kk[\Delta^t]\to\kk[\Delta^s]$. For $s\le t$ the
\emph{persistent Bass numbers} are
\[
\mu^{p,s\to t}_{U,G}=\rank\bigl(\mathbb M^p_{U,G}(f)_t\to\mathbb M^p_{U,G}(f)_s\bigr),
\]
\[
\mu^{p,s\to t}(\pp_U)=\rank\Bigl(\Ext^p_{S_{\pp_U}}(\kappa,\kk[\Delta^t]_{\pp_U})\to\Ext^p_{S_{\pp_U}}(\kappa,\kk[\Delta^s]_{\pp_U})\Bigr),
\]
with $\kappa=\kappa(\pp_U)$. The superscript $s\to t$ names the inclusion $\Delta^s\subseteq\Delta^t$; the maps
go from $t$ to $s$, as in \cite{SuwayyidWeiGraphs,HeSuwayyidWei}. (In \cite{HeSuwayyidWei} the letter $\lambda$
denoted persistent local cohomology numbers; here it denotes Lyubeznik numbers.) The persistent Bass numbers
at $\mathfrak q_U$ are those at $\pp_{U^c}$.
\end{definition}

\begin{theorem}\label{thm:bass-persistent}
Let $s\le t$, $d=p-|G|-1$ and $W=[n]\setminus(U\cup G)$. Then $\mathbb M^p_{U,G}(f)$ is isomorphic, as an
$\RR^{\mathrm{op}}$-indexed module, to $\bigl(\rH^d((\link_{\Delta^t}U)_W;\kk)\bigr)_t$ with restriction maps, and
\[
\mu^{p,s\to t}_{U,G}=\rank\Bigl(\rH^d\bigl((\link_{\Delta^t}U)_W;\kk\bigr)\to\rH^d\bigl((\link_{\Delta^s}U)_W;\kk\bigr)\Bigr),
\qquad
\mu^{p,s\to t}(\pp_U)=\sum_{G\subseteq U^c}\mu^{p,s\to t}_{U,G}.
\]
A void link contributes $0$.
\end{theorem}

\begin{proof}
By Lemma~\ref{lem:star-loc}(b),(c), the maps $\kk[\Delta^t]_{(U)}\to\kk[\Delta^s]_{(U)}$ are the quotients
$\kk[\link_{\Delta^t}U]\to\kk[\link_{\Delta^s}U]$ of Stanley--Reisner rings on $U^c$, with
$\link_{\Delta^s}U\subseteq\link_{\Delta^t}U$. Theorem~\ref{thm:bass-static}(b) with $V=U^c$ identifies the
degree-$(-\one_G)$ maps with restrictions, and the identifications commute with the maps, so they form an
isomorphism of persistence modules. For the total rank, the graded map
$\Ext^p(\kk,\kk[\Delta^t]_{(U)})\to\Ext^p(\kk,\kk[\Delta^s]_{(U)})$ is the direct sum of its graded pieces,
which vanish outside the degrees $-\one_G$ by Theorem~\ref{thm:bass-static}(a). Proposition~\ref{prop:bass-reduction}
says the ungraded map has the same rank.
\end{proof}

\begin{corollary}\label{cor:bass-betti}
\begin{enumerate}[label=\textup{(\alph*)},nosep]
\item For $U=\varnothing$ and $W=[n]\setminus G$,
$\mu^{p,s\to t}_{\varnothing,G}=\rank\bigl(\rH^{p-|G|-1}(\Delta^t_W)\to\rH^{p-|G|-1}(\Delta^s_W)\bigr)$.
This is the persistent Hochster formula of \cite{SuwayyidWeiGraphs} for the persistent graded Betti number in
homological degree $n-p$ and multidegree $\one_W$.
\item For general $U$, $\mu^{p,s\to t}_{U,G}$ is the persistent graded Betti number of the filtered link
$(\link_{\Delta^t}U)_t$, a filtration of $\link_\Delta U$ by the monotone function $\tau\mapsto f(\tau\cup U)$,
in homological degree $n-|U|-p$ and multidegree $\one_W$.
\end{enumerate}
\end{corollary}

\begin{proof}
(a) is the case $U=\varnothing$ of Theorem~\ref{thm:bass-persistent}. For (b), $\tau\in\link_{\Delta^t}U$ if and
only if $f(\tau\cup U)\le t$, so the filtered link is the sublevel filtration of $\tau\mapsto f(\tau\cup U)$, and
Theorem~\ref{thm:bass-persistent} together with the persistent Hochster formula of \cite{SuwayyidWeiGraphs}
applied to this filtration on the vertex set $U^c$ gives the identification.
\end{proof}

So at the maximal ideal nothing new appears, and at $\pp_U$ the persistent Bass numbers are the persistent
graded Betti numbers of the filtered link of $U$. The content of the construction is the choice of $U$ and the
localization, not a new algebraic invariant of the link.

\begin{theorem}\label{thm:bass-stability}
Each $\mathbb M^p_{U,G}(f)$ is of finite type with critical values in $\{f(\sigma):\sigma\in\Delta\}$, and so
has a finite barcode $\mathcal B^p_{U,G}(f)$. For monotone $f,g$,
$d_b\bigl(\mathcal B^p_{U,G}(f),\mathcal B^p_{U,G}(g)\bigr)\le\|f-g\|_\infty$. The same holds for the ungraded
modules at $\pp_U$, whose barcode is the union of the $\mathcal B^p_{U,G}(f)$ over $G\subseteq U^c$.
\end{theorem}

\begin{proof}
The assignment $\Delta'\mapsto\Ext^p_{S_{U^c}}(\kk,\kk[\Delta']_{(U)})_{-\one_G}$ is a contravariant functor on
subcomplexes of $\Delta$ with finite-dimensional values by Theorem~\ref{thm:bass-static}. Apply
Lemma~\ref{lem:functor-stability}. The ungraded module is the base change to $\kappa(\pp_U)$ of the direct sum
over $G$ by Proposition~\ref{prop:bass-reduction}, which gives the last statement.
\end{proof}

\begin{proposition}\label{prop:bass-bars}
Let $U\subseteq[n]$, $G\subseteq U^c$, $W=[n]\setminus(U\cup G)$, and $p=|G|$. If $U\in\Delta$, the barcode
$\mathcal B^{|G|}_{U,G}(f)$ is the single interval
\[
\Bigl[\,f(U),\ \min\{f(U\cup\{v\}):v\in W,\ U\cup\{v\}\in\Delta\}\Bigr),
\]
with $\min\varnothing=\infty$ and empty intervals omitted. If $U\notin\Delta$ the barcode is empty. In
particular, for $G=\varnothing$ and $p=0$ the barcode of $\mu^0(\pp_U,\kk[\Delta^\bullet])$ is the interval of
$t$ on which $U$ is a facet of $\Delta^t$. Moreover every bar of every $\mathcal B^p_{U,G}(f)$ is contained in
$[f(U),\infty)$.
\end{proposition}

\begin{proof}
Here $d=-1$. For a complex $K$, $\rH^{-1}(K)\neq0$ exactly when $K=\{\varnothing\}$, and then it is $\kk$.
The complex $(\link_{\Delta^t}U)_W$ is $\{\varnothing\}$ exactly when $U\in\Delta^t$ and no $v\in W$ has
$U\cup\{v\}\in\Delta^t$, that is, for $t$ in the displayed interval. Between two such values of $t$ the
restriction map is the identity of $\kk$. For $G=\varnothing$, $W=U^c$ and the condition says that $U$ is a
facet of $\Delta^t$; Corollary~\ref{cor:bass-static}(ii) shows that $\mu^0(\pp_U)=\mu^0_{U,\varnothing}$. The
last statement holds because $\link_{\Delta^t}U$ is void for $t<f(U)$.
\end{proof}

Since the associated primes of $\kk[\Delta^t]$ are the $\pp_U$ with $U$ a facet of $\Delta^t$, the degree-zero
Bass barcodes record the lifetime of each associated prime along the filtration. In degree
$p=|G|+1$ the modules are reduced $\rH^0$ of induced subcomplexes of links; for $G=\varnothing$ and
$U=\{i\}$ they count the connected components of $\link_{\Delta^t}(i)$ minus one whenever the link has a
vertex, and their bars can be read from a minimum spanning forest as for the degree-$0$ barcodes of
\cite{HeSuwayyidWei}. For $\mathfrak q_U$ the same statements hold with $U^c$ in place of $U$: the bars start at
$f(U^c)$ at the earliest, the degree-zero bar is the facet lifetime of $U^c$, and the degree-$1$ bars for
$G=\varnothing$ record the components of $\link_{\Delta^t}U^c$, a graph on part of $U$.

\section{Lyubeznik numbers as cohomology of a complex}\label{sec:lyub-static}

\subsection{Straight modules}

\begin{definition}[following Yanagawa \cite{Yanagawa2001}]
A $\ZZ^n$-graded $S$-module $N$ is \emph{straight} if $\dim_\kk N_a<\infty$ for all $a$ and, for every
$a\in\ZZ^n$ and every $j$ with $a_j\neq-1$, multiplication $x_j\colon N_a\to N_{a+\mathbf e_j}$ is bijective.
\end{definition}

\begin{lemma}\label{lem:straight}
\begin{enumerate}[label=\textup{(\alph*)},nosep]
\item For $U\subseteq[n]$, the localization $S_{x_U}$ is straight.
\item If $C^\bullet$ is a complex of straight modules with degree-preserving $S$-linear differentials, each
$H^i(C^\bullet)$ with finite-dimensional graded pieces is straight.
\end{enumerate}
\end{lemma}

\begin{proof}
(a) $(S_{x_U})_a$ is $\kk x^a$ if $a_j\ge0$ for all $j\notin U$, and $0$ otherwise. Let $a_j\neq-1$. If $j\in U$,
or $a_j\ge0$, both $(S_{x_U})_a$ and $(S_{x_U})_{a+\mathbf e_j}$ are $\kk x^a$ and $\kk x^{a+\mathbf e_j}$ or
both are zero, and $x_j$ matches them. If $j\notin U$ and $a_j\le-2$, both are zero.
(b) For fixed $a$ and $j$ with $a_j\ne-1$, multiplication by $x_j$ is an isomorphism of complexes of vector
spaces $C^\bullet_a\to C^\bullet_{a+\mathbf e_j}$, hence an isomorphism on cohomology.
\end{proof}

\begin{proposition}\label{prop:straight-ext}
Let $N$ be straight. Then $\Ext^p_S(\kk,N)_a=0$ for $a\neq-\one_{[n]}$, and
$\Ext^p_S(\kk,N)\cong H^p(B^\bullet(N))$, naturally in $N$, where
\[
B^p(N)=\bigoplus_{J\subseteq[n],\,|J|=p}N_{-\one_{[n]\setminus J}},\qquad
d(u_J)=\sum_{j\notin J}(-1)^{\nu_J(j)}\,x_ju_J\in N_{-\one_{[n]\setminus(J\cup\{j\})}}.
\]
\end{proposition}

\begin{proof}
If $a_j\neq-1$ for some $j$, then $(a+\one_J)_j=a_j\ne-1$ for $J\not\ni j$, and Lemma~\ref{lem:cone} shows
$K^\bullet(N)_a$ is acyclic. For $a=-\one_{[n]}$, $K^p(N)_a=\bigoplus_{|J|=p}N_{-\one_{[n]}+\one_J}e_J$, which is
$B^p(N)$ with the same differential.
\end{proof}

\subsection{The \v{C}ech complex on all nonfaces}

Let $\Delta\neq2^{[n]}$ and let $\Nf(\Delta)=\{N\subseteq[n]:N\notin\Delta\}$, a nonempty set. The monomials
$x_N$, $N\in\Nf(\Delta)$, generate $I_\Delta$ (if $\Delta$ is void, $x_\varnothing=1$). Fix once and for all a
total order on $2^{[n]}$. For $A\subseteq\Nf(\Delta)$ let $U(A)=\bigcup_{N\in A}N$. The \v{C}ech complex of
$S$ on these generators is
\[
\check C^i(\Delta)=\bigoplus_{A\subseteq\Nf(\Delta),\,|A|=i}S_{x_{U(A)}}\,e_A,\qquad
d(u\,e_A)=\sum_{N\in\Nf(\Delta)\setminus A}(-1)^{\nu_A(N)}u\,e_{A\cup\{N\}},
\]
where $u$ is mapped by the localization map. Its cohomology is $\Hloc{i}{I_\Delta}(S)$ \cite[\S3.5]{BrunsHerzog}.
By Lemma~\ref{lem:straight}, $\Hloc{i}{I_\Delta}(S)$ is straight; its graded pieces are finite-dimensional since
each term of $\check C^\bullet(\Delta)$ has graded pieces of dimension at most $1$ and there are finitely many
terms.

\begin{lemma}\label{lem:cech-maps}
Let $\Delta'\subseteq\Delta''$ be subcomplexes of $\Delta$, both different from $2^{[n]}$. Then
$\Nf(\Delta'')\subseteq\Nf(\Delta')$, and the projection $\pi\colon\check C^\bullet(\Delta')\to\check C^\bullet(\Delta'')$
onto the summands with $A\subseteq\Nf(\Delta'')$ is a chain map. On cohomology it induces the natural map
$\Hloc{i}{I_{\Delta'}}(S)\to\Hloc{i}{I_{\Delta''}}(S)$ coming from $\Gamma_{I_{\Delta'}}\subseteq\Gamma_{I_{\Delta''}}$.
Projections compose.
\end{lemma}

\begin{proof}
Write $B'=\Nf(\Delta')\supseteq B''=\Nf(\Delta'')$. If $A\subseteq B''$, then $\pi d(u\,e_A)$ keeps the terms with
$N\in B''$, which is $d\pi(u\,e_A)$. If $A\not\subseteq B''$, both sides are $0$. The same formulas define
$\pi$ for $\check C^\bullet(B';M)\to\check C^\bullet(B'';M)$ for every $S$-module $M$, naturally in $M$. Since
localization is exact, both $M\mapsto H^i(\check C^\bullet(B';M))$ and $M\mapsto H^i(\check C^\bullet(B'';M))$ are
$\delta$-functors, isomorphic to $\Hloc{i}{I_{\Delta'}}$ and $\Hloc{i}{I_{\Delta''}}$ \cite[\S3.5]{BrunsHerzog},
and $\pi$ commutes with connecting maps, so it is a morphism of $\delta$-functors. In degree $0$,
$H^0(\check C^\bullet(B;M))=\Gamma_{(x_N:N\in B)}(M)$ and $\pi$ is the identity of $M$, so it induces the
inclusion $\Gamma_{I_{\Delta'}}(M)\subseteq\Gamma_{I_{\Delta''}}(M)$. Local cohomology is a universal
$\delta$-functor, so a morphism of $\delta$-functors out of it is determined by its degree-$0$ part
\cite[\S2.1]{Weibel}.
\end{proof}

\subsection{Graded pieces of local cohomology}

\begin{proposition}\label{prop:graded-pieces}
Let $\Delta\neq2^{[n]}$ and $N^i=\Hloc{i}{I_\Delta}(S)$.
\begin{enumerate}[label=\textup{(\alph*)},nosep]
\item $N^i_{\mathbf 0}=0$.
\item For $\varnothing\neq G\subseteq[n]$ there is an isomorphism
$\theta_G\colon N^i_{-\one_G}\to\rH^{\,i-2}\bigl((\Delta^\vee)_G;\kk\bigr)$.
\item For $j\in G$ with $G\ne\{j\}$, $\theta_{G\setminus\{j\}}\circ x_j=\res\circ\theta_G$, where $\res$ is
restriction along $(\Delta^\vee)_{G\setminus\{j\}}\subseteq(\Delta^\vee)_G$.
\item For $\Delta'\subseteq\Delta''$ as in Lemma~\ref{lem:cech-maps}, the natural map
$\Hloc{i}{I_{\Delta'}}(S)_{-\one_G}\to\Hloc{i}{I_{\Delta''}}(S)_{-\one_G}$ corresponds under $\theta_G$ to restriction
along $(\Delta''^\vee)_G\subseteq(\Delta'^\vee)_G$.
\end{enumerate}
\end{proposition}

Parts (b) and (c) are due to Musta\c{t}\u{a} \cite[Cor.~2.2(a),(b)]{Mustata}, whose
complex ``$\Delta$'' is our $\Delta^\vee$. We prove them again because the proof below
makes the identifications natural in $\Delta$, which gives the new part (d).

\begin{proof} 
\emph{Step 1: the degree-$(-\one_G)$ part of the \v{C}ech complex.} Since $(S_{x_U})_{-\one_G}$ is $\kk x^{-\one_G}$
if $G\subseteq U$ and $0$ otherwise, and localization maps send $x^{-\one_G}$ to itself,
$\check C^i(\Delta)_{-\one_G}$ has basis $\{e_A:|A|=i,\ G\subseteq U(A)\}$ and $d(e_A)=\sum_N(-1)^{\nu_A(N)}e_{A\cup\{N\}}$.
Let $Y=2^{\Nf(\Delta)}$ be the full simplex on $\Nf(\Delta)$ and
\[
X_G=\{A\subseteq\Nf(\Delta):G\not\subseteq U(A)\},
\]
a subcomplex. Sending $e_A$ to the cochain $A^\vee$ identifies $\check C^i(\Delta)_{-\one_G}$ with
$\rC^{\,i-1}(Y,X_G)$, the reduced cochains of $Y$ vanishing on $X_G$; the signs agree with the coboundary. Since
$\Nf(\Delta)\neq\varnothing$, $\rC^\bullet(Y)$ is acyclic.

(a) For $G=\varnothing$, $X_G$ is void, so the complex is $\rC^{\bullet-1}(Y)$, which is acyclic.

\emph{Step 2: connecting map.} For $G\neq\varnothing$, $\varnothing\in X_G$. The exact sequence
$0\to\rC^\bullet(Y,X_G)\to\rC^\bullet(Y)\to\rC^\bullet(X_G)\to0$ and acyclicity of $\rC^\bullet(Y)$ give an
isomorphism $\partial\colon\rH^{i-2}(X_G)\to H^{i-1}(\rC^\bullet(Y,X_G))\cong N^i_{-\one_G}$, natural for maps of pairs.

\emph{Step 3: nerve.} For $g\in G$ let $\Nf_g=\{N\in\Nf(\Delta):g\notin N\}$ and $Y_g=2^{\Nf_g}$. A set $A$ lies in
$X_G$ exactly when $A\subseteq\Nf_g$ for some $g\in G$, so $X_G=\bigcup_{g\in G}Y_g$. For $\varnothing\neq H\subseteq G$,
$\bigcap_{g\in H}Y_g$ is the full simplex on $\{N\in\Nf(\Delta):N\cap H=\varnothing\}$. This set is nonempty exactly
when $[n]\setminus H\notin\Delta$, that is, when $H\in\Delta^\vee$. So the geometric realizations $|Y_g|$ that are
nonempty form a closed cover of $|X_G|$ with contractible nonempty intersections, and its nerve is
$(\Delta^\vee)_G$. By the nerve theorem \cite{Bjorner}, $|X_G|\simeq|(\Delta^\vee)_G|$ when these are nonempty. They
are empty at the same time, namely when no $\Nf_g$ is nonempty, and then $X_G=\{\varnothing\}=(\Delta^\vee)_G$
(recall $\varnothing\in\Delta^\vee$ since $[n]\notin\Delta$). This gives $\rH^{i-2}(X_G)\cong\rH^{i-2}((\Delta^\vee)_G)$;
composing with $\partial^{-1}$ defines $\theta_G$ and proves (b).

\emph{Step 4: maps.} (c) Multiplication by $x_j$ sends $x^{-\one_G}$ to $x^{-\one_{G\setminus\{j\}}}$ in every
$S_{x_U}$ with $G\subseteq U$, so on \v{C}ech pieces it sends $e_A$ to $e_A$. This is the inclusion
$\rC^\bullet(Y,X_G)\subseteq\rC^\bullet(Y,X_{G\setminus\{j\}})$ induced by the map of pairs
$(Y,X_{G\setminus\{j\}})\to(Y,X_G)$. By naturality of $\partial$ it corresponds to restriction
$\rH(X_G)\to\rH(X_{G\setminus\{j\}})$. The cover of $X_{G\setminus\{j\}}$ is the subfamily of the cover of $X_G$ indexed
by $G\setminus\{j\}$, and the nerve inclusion is $(\Delta^\vee)_{G\setminus\{j\}}\subseteq(\Delta^\vee)_G$. The
functorial nerve theorem \cite{BauerKerberRollRolle} says the nerve equivalences commute up to homotopy with these
inclusions, which gives (c) on cohomology.

(d) Let $Y',X'_G$ and $Y'',X''_G$ be the complexes above for $\Delta'$ and $\Delta''$. Then $Y''\subseteq Y'$,
$X''_G=X'_G\cap Y''$, and $Y''_g=Y'_g\cap Y''$. By Lemma~\ref{lem:cech-maps} the natural map is induced by $\pi$,
which sends $e_A$ to $e_A$ if $A\subseteq\Nf(\Delta'')$ and to $0$ otherwise. This is restriction of cochains along
the map of pairs $(Y'',X''_G)\subseteq(Y',X'_G)$. Naturality of $\partial$ turns it into restriction
$\rH(X'_G)\to\rH(X''_G)$. The covers satisfy $Y''_g\subseteq Y'_g$ with the same index set, and the nerves are
$(\Delta''^\vee)_G\subseteq(\Delta'^\vee)_G$, so the functorial nerve theorem gives (d). When the realizations are
empty the statements are checked directly, since all complexes involved are $\{\varnothing\}$ or the relevant groups
vanish.
\end{proof}

\subsection{The Lyubeznik complex}
An expert in the domain might also recognize the complex as dual of the Àlvarez Montaner–Vahidi defined and discussed in \cite[\S4]{AlvarezMontanerVahidi} via $r = n - q$.
\begin{definition}\label{def:lyub-complex}
Let $\Delta\neq2^{[n]}$ and $0\le q\le n$. The \emph{Lyubeznik complex} $\Lc^\bullet_q(\Delta)$ is
\[
\Lc^p_q(\Delta)=\bigoplus_{\substack{J\subseteq[n],\ |J|=p\\ J\neq[n]}}\rH^{\,n-q-2}\bigl((\Delta^\vee)_{[n]\setminus J};\kk\bigr),
\qquad
d\bigl(u_J\bigr)=\sum_{\substack{j\notin J\\ J\cup\{j\}\neq[n]}}(-1)^{\nu_J(j)}\,u_J\big|_{(\Delta^\vee)_{[n]\setminus(J\cup\{j\})}} .
\]
Up to $\kk$-duality and notation, the complex defined next is the complex of \`Alvarez
Montaner and Vahidi \cite[Prop.~3.3, Cor.~3.4 and \S4]{AlvarezMontanerVahidi} for the
module $\Hloc{r}{I_\Delta}(S)$ with $r=n-q$, and the equality
$\lambda_{p,q}(\kk[\Delta])=\dim_\kk H^p(\Lc^\bullet_q(\Delta))$ is known. What is new in
Theorem~\ref{thm:lyub-complex} is the naturality in $\Delta$.
\end{definition}

By Lemma~\ref{lem:dual-of-link}, the summand of $J$ is zero unless $J\in\Delta$, and it depends only on
$\link_\Delta J$: it is $\rH^{n-q-2}$ of the Alexander dual of $\link_\Delta J$ in $[n]\setminus J$.

\begin{theorem}\label{thm:lyub-complex}
Let $\Delta\neq2^{[n]}$. For all $p,q$,
\[
\lambda_{p,q}(\kk[\Delta])=\dim_\kk H^p\bigl(\Lc^\bullet_q(\Delta)\bigr).
\]
More precisely, there is an isomorphism $\Ext^p_S(\kk,\Hloc{n-q}{I_\Delta}(S))\cong H^p(\Lc^\bullet_q(\Delta))$, and for
$\Delta'\subseteq\Delta''$ as above it carries the map induced by $\Hloc{n-q}{I_{\Delta'}}(S)\to\Hloc{n-q}{I_{\Delta''}}(S)$
to the map induced by the chain map of Definition~\ref{def:lyub-complex}. For $\Delta=2^{[n]}$ one has
$\lambda_{n,n}=1$ and all other Lyubeznik numbers vanish.
\end{theorem}

\begin{proof}
Let $N=\Hloc{n-q}{I_\Delta}(S)$, which is straight. By Section~\ref{sec:background}, $\lambda_{p,q}=\dim\Ext^p_S(\kk,N)$,
and by Proposition~\ref{prop:straight-ext} this is $\dim H^p(B^\bullet(N))$. The summand of $B^p(N)$ for $J=[n]$ is
$N_{\mathbf 0}=0$ by Proposition~\ref{prop:graded-pieces}(a). For $J\neq[n]$ apply $\theta_{[n]\setminus J}$ with
$i=n-q$, so $i-2=n-q-2$. By Proposition~\ref{prop:graded-pieces}(c) these isomorphisms turn $x_j$ into restriction;
the signs $(-1)^{\nu_J(j)}$ are the same on both sides, and terms with $J\cup\{j\}=[n]$ vanish. So
$\theta$ is an isomorphism of complexes $B^\bullet(N)\to\Lc^\bullet_q(\Delta)$. Naturality in $\Delta$ follows from
Proposition~\ref{prop:graded-pieces}(d) and the naturality in Proposition~\ref{prop:straight-ext}. For
$\Delta=2^{[n]}$, $I_\Delta=0$, $\Hloc{0}{0}(S)=S$ and the other local cohomology vanishes; $S$ is straight and
$B^p(S)$ is nonzero only for $J=[n]$, where it is $S_{\mathbf 0}=\kk$, giving $\lambda_{n,n}=1$.
\end{proof}

\begin{corollary}\label{cor:lyub-terms}
Let $\Delta\ne2^{[n]}$, $d=\dim\kk[\Delta]$ (the largest size of a face), and $0\le p,q\le n$.
\begin{enumerate}[label=\textup{(\alph*)},nosep]
\item $\dim_\kk\Lc^p_q(\Delta)=\sum_{J\in\Delta,\ |J|=p}\dim_\kk\rH_{q-p-1}(\link_\Delta J;\kk)$.
\item $\lambda_{p,q}(\kk[\Delta])=0$ unless $p\le q\le d$.
\item $\sum_p(-1)^p\lambda_{p,q}(\kk[\Delta])=\sum_{J\in\Delta}(-1)^{|J|}\dim_\kk\rH_{q-|J|-1}(\link_\Delta J;\kk)$.
\item $\lambda_{q,q}(\kk[\Delta])$ is the dimension of the cokernel of $d\colon\Lc^{q-1}_q\to\Lc^q_q$, and
$\Lc^q_q(\Delta)$ has one copy of $\kk$ for each facet of $\Delta$ with $q$ vertices. In particular
$\lambda_{q,q}(\kk[\Delta])$ is at most the number of such facets.
\item For $n \geq 2$, $\lambda_{0,q}(\kk[\Delta])$ is the dimension of the common kernel of the restriction maps
$\rH^{n-q-2}(\Delta^\vee)\to\rH^{n-q-2}((\Delta^\vee)_{[n]\setminus\{j\}})$, $j\in[n]$; here
$\dim\rH^{n-q-2}(\Delta^\vee)=\dim\rH_{q-1}(\Delta)$.
\end{enumerate}
\end{corollary}

\begin{proof}
(a) Let $J\in\Delta$, $J\ne[n]$, and $V=[n]\setminus J$. By Lemma~\ref{lem:dual-of-link}, $(\Delta^\vee)_V$ is the
Alexander dual of $\link_\Delta J$ in $V$, and it is not $2^V$ since $J\in\Delta$. Theorem~\ref{thm:alexander} with
$i=n-q-2$ and $|V|=n-p$ gives $\rH^{n-q-2}((\Delta^\vee)_V)\cong\rH_{q-p-1}(\link_\Delta J)$. For $J\notin\Delta$ the
summand is the cohomology of a full simplex on $V\ne\varnothing$, which is $0$.
(b) A nonzero term needs $-1\le q-p-1\le\dim\link_\Delta J\le d-p-1$.
(c) Take the Euler characteristic of $\Lc^\bullet_q$ using (a).
(d) $\Lc^{q+1}_q=0$ since it involves $\rH_{-2}$. A summand of $\Lc^q_q$ is $\rH_{-1}(\link_\Delta J)$, which is $\kk$
exactly when $J$ is a facet.
(e) $\Lc^{-1}_q=0$, and $\Lc^0_q=\rH^{n-q-2}(\Delta^\vee)$ with differential the sum of the restrictions (all signs
are $+1$). The last statement is Theorem~\ref{thm:alexander} for $J=\varnothing$.
\end{proof}

Corollary~\ref{cor:lyub-terms}(a) describes the terms through links. Section~\ref{sec:link-model} gives an
equivalent costar model and describes the differential in terms of links, with signs: the map from the
summand indexed by $J$ to that indexed by $J\cup\{j\}$ deletes the vertex $j$ from a face of
$\link_\Delta J$ and moves it into $J$ (Lemma~\ref{lem:costar-link}).

\begin{example}\label{ex:two-planes}
Let $\Delta$ have facets $\{1,2\}$ and $\{3,4\}$, so $\kk[\Delta]$ is the coordinate ring of two planes in
$\mathbb A^4$ meeting at the origin. For $q=1$ the only nonzero term is $\Lc^0_1\cong\rH_0(\Delta)=\kk$, so
$\lambda_{0,1}=1$. For $q=2$ the only nonzero term is $\Lc^2_2=\kk^2$, one copy for each facet, so
$\lambda_{2,2}=2$. The other terms vanish because every vertex link is a point and every edge link is
$\{\varnothing\}$. This is the known Lyubeznik table of this ring.
\end{example}

\begin{example}\label{ex:rp2}
Let $\Delta$ be the six-vertex triangulation of the real projective plane, with facets
$012,023,034,045,051,124,235,341,452,513$ on $\{0,\dots,5\}$. The link of each vertex is a $5$-cycle and the link of
each edge is two points. If $\operatorname{char}\kk\ne2$, then $\rH_1(\Delta)=\rH_2(\Delta)=0$, so $\Lc^0_q=0$ for all
$q$, and $\kk[\Delta]$ is Cohen--Macaulay with $\lambda_{3,3}=1$ and all other entries zero. If $\operatorname{char}\kk=2$,
then $\Lc^0_2=\rH_1(\Delta)=\kk$ while $\Lc^1_2=\bigoplus_v\rH_0(\link v)=0$, so $\lambda_{0,2}=1$. Computing the
rest of the complex gives $\lambda_{0,2}=\lambda_{2,3}=\lambda_{3,3}=1$ and zero otherwise. For $q=3$ the term
dimensions are $1,6,15,10$, consistent with Corollary~\ref{cor:lyub-terms}(c). So the Lyubeznik table depends on the
characteristic. For characteristic dependence of Lyubeznik numbers of Stanley--Reisner rings in general, see
\cite{AlvarezMontanerVahidi}.
\end{example}

\section{Persistent Lyubeznik numbers}\label{sec:lyub-persistent}

\begin{definition}\label{def:pers-lyub}
Let $f$ be monotone. The \emph{Lyubeznik persistence module} is the $\RR$-indexed module
\[
\boldsymbol\Lambda^{p,q}(f)=\Bigl(\Ext^p_S\bigl(\kk,\Hloc{n-q}{I_{\Delta^t}}(S)\bigr)\Bigr)_{t\in\RR},
\]
with maps from $s$ to $t$ ($s\le t$) induced by $\Hloc{n-q}{I_{\Delta^s}}(S)\to\Hloc{n-q}{I_{\Delta^t}}(S)$. The
\emph{persistent Lyubeznik numbers} are $\lambda^{s\to t}_{p,q}=\rank\bigl(\boldsymbol\Lambda^{p,q}(f)_s\to\boldsymbol\Lambda^{p,q}(f)_t\bigr)$.
For $s=t$ this is $\lambda_{p,q}(\kk[\Delta^t])$. Unlike the Bass modules, these maps go in the direction of the
filtration.
\end{definition}

For $t<f(\varnothing)$ we have $I_{\Delta^t}=S$ and the module is $0$. If $\Delta=2^{[n]}$, the top level
$\Delta^t=2^{[n]}$ has $\lambda_{n,n}=1$ and nothing else, while every other level has $\lambda_{n,n}=0$ by
Corollary~\ref{cor:lyub-terms}(b); so all maps between the top level and a lower level have rank $0$. We therefore
state results for levels different from $2^{[n]}$.

\begin{theorem}\label{thm:lyub-persistent}
Let $s\le t$ with $\Delta^s,\Delta^t\neq2^{[n]}$. Then
\[
\lambda^{s\to t}_{p,q}=\rank\Bigl(H^p\bigl(\Lc^\bullet_q(\Delta^s)\bigr)\to H^p\bigl(\Lc^\bullet_q(\Delta^t)\bigr)\Bigr),
\]
the map induced by restriction along $(\Delta^t)^\vee_{[n]\setminus J}\subseteq(\Delta^s)^\vee_{[n]\setminus J}$ in each
summand. Equivalently, $\boldsymbol\Lambda^{p,q}(f)$ is isomorphic, on levels different from $2^{[n]}$, to the $p$-th
cohomology of the complex of persistence modules $t\mapsto\Lc^\bullet_q(\Delta^t)$.
\end{theorem}

\begin{proof}
Theorem~\ref{thm:lyub-complex} applied to $\Delta'=\Delta^s\subseteq\Delta''=\Delta^t$.
\end{proof}

\begin{theorem}\label{thm:lyub-stability}
Each $\boldsymbol\Lambda^{p,q}(f)$ is of finite type with critical values in $\{f(\sigma):\sigma\in\Delta\}$, so it has
a finite barcode $\mathcal B_{p,q}(f)$, and $d_b\bigl(\mathcal B_{p,q}(f),\mathcal B_{p,q}(g)\bigr)\le\|f-g\|_\infty$.
\end{theorem}

\begin{proof}
$\Delta'\mapsto\Ext^p_S(\kk,\Hloc{n-q}{I_{\Delta'}}(S))$ is a covariant functor on subcomplexes of $\Delta$: the maps come
from inclusions of torsion functors, which compose. Its values are finite-dimensional by
Theorem~\ref{thm:lyub-complex}. Apply Lemma~\ref{lem:functor-stability}.
\end{proof}

\begin{proposition}\label{prop:top-row}
Let $s\le t$ with $\Delta^s,\Delta^t\ne2^{[n]}$. Then
\[
\lambda_{q,q}(\kk[\Delta^t])-\lambda^{s\to t}_{q,q}\le\bigl|\{J:\ |J|=q,\ s<f(J)\le t,\ J\text{ a facet of }\Delta^t\}\bigr|.
\]
Consequently, every bar of $\mathcal B_{q,q}(f)$ has left endpoint $f(J)$ for some face $J$ with $q$ vertices.
\end{proposition}

\begin{proof}
Write $V_u=H^q(\Lc^\bullet_q(\Delta^u))$, which is the cokernel of $d\colon\Lc^{q-1}_q(\Delta^u)\to\Lc^q_q(\Delta^u)$ by
Corollary~\ref{cor:lyub-terms}(d). The map $V_s\to V_t$ is induced by the chain map $\varphi\colon\Lc^q_q(\Delta^s)\to\Lc^q_q(\Delta^t)$,
so $\dim V_t-\rank(V_s\to V_t)\le\dim\Lc^q_q(\Delta^t)-\dim\varphi(\Lc^q_q(\Delta^s))$. The summand of $J$ in
$\Lc^q_q(\Delta^t)$ is nonzero only if $J$ is a facet of $\Delta^t$. If $J$ is a facet of both $\Delta^s$ and $\Delta^t$,
then $\link_{\Delta^s}J=\link_{\Delta^t}J=\{\varnothing\}$, so by Lemma~\ref{lem:dual-of-link} the complexes
$(\Delta^s)^\vee_{[n]\setminus J}$ and $(\Delta^t)^\vee_{[n]\setminus J}$ are equal and $\varphi$ is the identity on this
summand. If $J$ is a facet of $\Delta^t$ but not of $\Delta^s$, then $J\notin\Delta^s$, since a coface of $J$ in $\Delta^s$
would lie in $\Delta^t$; so $s<f(J)\le t$. This proves the inequality.

The left side counts the bars that contain $t$ and not $s$. If a bar has left endpoint $a$, choose $s<a$ such that
$f$ takes no value in $[s,a)$ and take $t=a$. The right side counts faces with $f(J)=a$ and $q$ vertices, so there is
at least one.
\end{proof}

\section{Costar and face-link complexes}\label{sec:link-model}

The description of multigraded local cohomology of
Stanley--Reisner rings through the relative cohomology
of costars, together with the corresponding multiplication
maps, goes back to Gr\"abe \cite{Grabe1984}.
The corresponding link--costar homology isomorphism
also appears in his work. The complex $\Lc_q^\bullet(\Delta)$ is defined using Alexander duals. We give an equivalent description in terms of relative simplicial chains of $\Delta$ and its costars. Its differentials are induced by quotient maps, which makes their behavior under inclusions of complexes immediate and permits computation without forming Alexander duals or minimal free resolutions. The relative groups are also the reduced homology groups of links, so the construction describes how the local homology at different faces fits together. This distinction matters for persistence: a filtration may change the links of faces without changing the reduced homology of the ambient complex.

\subsection{Relative homology of costars}

For $F\subseteq[n]$ put
\[
Q_\bullet(F;\Delta)=\rC_\bullet(\Delta;\kk)/\rC_\bullet(\cost_\Delta F;\kk),\qquad
\Kc_q(F;\Delta)=H_{q-1}(Q_\bullet(F;\Delta)).
\]
Here $\rC_\bullet$ is the augmented chain complex, whereas $\rC^\bullet$ denotes reduced cochains. In particular, $Q_\bullet(\varnothing;\Delta)=\rC_\bullet(\Delta;\kk)$.

For $F\subseteq G\subseteq[n]$ we have $\cost_\Delta F\subseteq\cost_\Delta G$, so the identity on $\rC_\bullet(\Delta;\kk)$ induces a quotient chain map
\[
P_{F,G}^\Delta:Q_\bullet(F;\Delta)\longrightarrow Q_\bullet(G;\Delta),\qquad
P_{F,G}^\Delta([S])=\begin{cases}[S],&G\subseteq S,\\0,&G\not\subseteq S,\end{cases}
\]
for a face $S$ containing $F$. Let $\rho_{F,G}^\Delta=H_{q-1}(P_{F,G}^\Delta):\Kc_q(F;\Delta)\to\Kc_q(G;\Delta)$. Since $P_{G,H}^\Delta P_{F,G}^\Delta=P_{F,H}^\Delta$, the induced maps satisfy $\rho_{G,H}^\Delta\rho_{F,G}^\Delta=\rho_{F,H}^\Delta$. They also commute with the maps induced by inclusions $\Delta'\subseteq\Delta''$, as is clear from the corresponding inclusions of pairs.

\begin{definition}\label{def:face-link-complex}
For $0\le q\le n$ put
\[
\Kc_q^p(\Delta)=\bigoplus_{\substack{F\subseteq[n]\\|F|=p}}\Kc_q(F;\Delta).
\]
The differential of $\Kc_q^\bullet(\Delta)$ is given, for $|G|=p+1$, by
\[
(dz)_G=\sum_{j\in G}(-1)^{\nu_{G\setminus\{j\}}(j)}
\rho^\Delta_{G\setminus\{j\},G}(z_{G\setminus\{j\}}),
\]
where $\nu_F(j)=|\{u\in F:u<j\}|$. A nonface $F$ contributes zero.
\end{definition}

\begin{lemma}\label{lem:face-link-differential}
The differential of Definition~\ref{def:face-link-complex} satisfies $d^2=0$. For $\Delta'\subseteq\Delta''$, the inclusions of pairs $(\Delta',\cost_{\Delta'}F)\hookrightarrow(\Delta'',\cost_{\Delta''}F)$ induce a cochain map $\Kc_q^\bullet(\Delta')\to\Kc_q^\bullet(\Delta'')$, and these maps respect composition.
\end{lemma}

\begin{proof}
Write $D_p^\Delta=(d_{G,F}^\Delta)_{|G|=p+1,\,|F|=p}$ for the block matrix of the differential, where
\[
d_{G,F}^\Delta=\begin{cases}
(-1)^{\nu_F(j)}\rho^\Delta_{F,G},&G=F\cup\{j\},\\
0,&\text{otherwise}.
\end{cases}
\]
The $(H,F)$-block of $D_{p+1}^\Delta D_p^\Delta$ is $\sum_{|G|=p+1}d_{H,G}^\Delta d_{G,F}^\Delta$, and is zero unless $F\subset H$. In the latter case write $H\setminus F=\{i,j\}$ with $i<j$ and put $G_i=F\cup\{i\}$, $G_j=F\cup\{j\}$. The block equals
\[
\begin{aligned}
&(-1)^{\nu_F(i)+\nu_{G_i}(j)}\rho^\Delta_{G_i,H}\rho^\Delta_{F,G_i}
+(-1)^{\nu_F(j)+\nu_{G_j}(i)}\rho^\Delta_{G_j,H}\rho^\Delta_{F,G_j}.
\end{aligned}
\]
Both compositions are $\rho^\Delta_{F,H}$, whereas $\nu_{G_i}(j)=\nu_F(j)+1$ and $\nu_{G_j}(i)=\nu_F(i)$. The terms cancel, so $D_{p+1}^\Delta D_p^\Delta=0$.

For $\Delta'\subseteq\Delta''$ let $a_F:\Kc_q(F;\Delta')\to\Kc_q(F;\Delta'')$ be induced by the inclusion of pairs, and set $A_p=\operatorname{diag}(a_F)_{|F|=p}$. Naturality of relative homology gives
\[
a_G\rho^{\Delta'}_{F,G}=\rho^{\Delta''}_{F,G}a_F.
\]
Comparing blocks yields $A_{p+1}D_p^{\Delta'}=D_p^{\Delta''}A_p$, so $A_\bullet$ is a cochain map. The identity $A_p^{\Delta',\Delta'''}=A_p^{\Delta'',\Delta'''}A_p^{\Delta',\Delta''}$ follows from composition of inclusions of pairs.
\end{proof}

The first isomorphism in the following lemma is the classical
link--costar isomorphism of Gr\"abe
\cite[\S 1.3]{Grabe1984}. We also describe the maps between
the face-indexed summands and their compatibility with
inclusions of simplicial complexes.

\begin{lemma}\label{lem:costar-link}
For $F\in\Delta$ and $p=|F|$, the correspondence $[T]\mapsto[F,T]$ induces an isomorphism
\[
\theta_F:\rH_{q-p-1}(\link_\Delta F;\kk)\xrightarrow{\ \cong\ }\Kc_q(F;\Delta).
\]
If $F\notin\Delta$, then $\Kc_q(F;\Delta)=0$. Put
\[
\mathcal H_q^p(\Delta)=\bigoplus_{\substack{F\in\Delta\\|F|=p}}\rH_{q-p-1}(\link_\Delta F;\kk).
\]
For $G=F\cup\{j\}$ define a degree-lowering map on augmented link chains by
\[
\tau_{F,j}([T])=\begin{cases}(-1)^{\nu_T(j)}[T\setminus\{j\}],&j\in T,\\0,&j\notin T.\end{cases}
\]
Then $\partial\tau_{F,j}=-\tau_{F,j}\partial$, and
\[
(d_{\mathcal H}z)_G=(-1)^p\sum_{j\in G}(\tau_{G\setminus\{j\},j})_*(z_{G\setminus\{j\}}),\qquad |G|=p+1,
\]
defines a differential on $\mathcal H_q^\bullet(\Delta)$. The maps $\theta_F$ give a natural isomorphism $\mathcal H_q^\bullet(\Delta)\cong\Kc_q^\bullet(\Delta)$ of cochain complexes.
\end{lemma}

\begin{proof}
Let $F\in\Delta$ and $p=|F|$. Every face containing $F$ is uniquely $F\cup T$ for $T\in\link_\Delta F$. In degree $r$ define
\[
\Theta_F:\rC_{r-p}(\link_\Delta F;\kk)\longrightarrow Q_r(F;\Delta),\qquad
\Theta_F([T])=(-1)^{p(r-q+1)}[F,T].
\]
This is a bijection of bases. Deleting a vertex of $F$ gives zero in the quotient, while deleting a vertex of $T$ gives its augmented boundary with sign $(-1)^p$. The normalization thus gives $\partial_Q\Theta_F=\Theta_F\partial_{\link}$. Hence $\Theta_F$ induces the stated isomorphism $\theta_F([z])=[\Theta_F(z)]$. In degree $r=q-1$ the normalization is $1$, so $\theta_F$ is induced by $[T]\mapsto[F,T]$. If $F\notin\Delta$, no face contains $F$ and $Q_\bullet(F;\Delta)=0$.

Let $G=F\cup\{j\}\in\Delta$. The map $\rho^\Delta_{F,G}$ kills $[F,T]$ when $j\notin T$ and retains it otherwise. If $j\in T$, moving $j$ to its increasing position in $G$ gives
\[
[F,T]=(-1)^{p-\nu_F(j)+\nu_T(j)}[G,T\setminus\{j\}].
\]
After multiplication by the incidence sign this becomes
\[
(-1)^{\nu_F(j)}[F,T]=(-1)^{p+\nu_T(j)}[G,T\setminus\{j\}].
\]
Thus $\theta_G^{-1}\bigl((-1)^{\nu_F(j)}\rho^\Delta_{F,G}\bigr)\theta_F=(-1)^p(\tau_{F,j})_*$. The identity $\partial\tau_{F,j}=-\tau_{F,j}\partial$ follows by applying the augmented simplicial boundary twice, and shows that $\tau_{F,j}$ induces a map on homology. Summing over $j\in G$ gives the stated differential on $\mathcal H_q^\bullet(\Delta)$ and the isomorphism with $\Kc_q^\bullet(\Delta)$. The identifications are natural under simplicial inclusions because the oriented face correspondences commute with the inclusions of relative chain complexes.
\end{proof}

\begin{theorem}\label{thm:link-model}
Let $\Delta\ne2^{[n]}$. For all $p,q$ there are natural isomorphisms
\[
H^p(\Lc_q^\bullet(\Delta))\cong H^p(\Kc_q^\bullet(\Delta)),\qquad
\lambda_{p,q}(\kk[\Delta])=\dim_\kk H^p(\Kc_q^\bullet(\Delta)).
\]
More precisely, $\Lc_q^\bullet(\Delta)$ and $\Kc_q^\bullet(\Delta)$ are isomorphic as cochain complexes, compatibly with the maps induced by $\Delta'\subseteq\Delta''$ when both complexes are proper. The face-link complex also computes the Lyubeznik numbers for $\Delta=2^{[n]}$.
\end{theorem}

\begin{proof}
Put $N=[n]$, $d=n-q-2$, and, for $F\subsetneq N$, let $W_F=N\setminus F$ and $D_F=(\Delta^\vee)_{W_F}$. Since $\Delta\ne2^N$, the complex $D_F$ contains $\varnothing$ and $W_F\ne\varnothing$. The exact sequence of augmented cochain complexes
\[
0\longrightarrow C^\bullet(2^{W_F},D_F)\longrightarrow\rC^\bullet(2^{W_F})\longrightarrow\rC^\bullet(D_F)\longrightarrow0
\]
and the acyclicity of $2^{W_F}$ give a natural connecting isomorphism
\[
\kappa_F:\rH^d(D_F;\kk)\xrightarrow{\ \cong\ }H^{d+1}(2^{W_F},D_F;\kk).
\]

We identify this relative cohomology with the costar relative homology. For $H\subseteq W_F$, we have $H\notin D_F$ if and only if $S=N\setminus H\in\Delta$, in which case $F\subseteq S$. On relative cochain bases put
\[
E_F^r:C^r(2^{W_F},D_F)\longrightarrow Q_{n-r-2}(F;\Delta),\qquad
E_F^r(e_H^*)=\varepsilon(H)[N\setminus H],\quad
\varepsilon(H)=(-1)^{\sum_{h\in H}(h-1)}.
\]
This is a bijection of bases satisfying $\partial E_F^r=E_F^{r+1}\delta$. Indeed, deleting a vertex of $F$ gives zero in the quotient, and for $v\in S\setminus F$ the remaining coefficients agree because
\[
(-1)^{\nu_H(v)}\varepsilon(H\cup\{v\})=(-1)^{\nu_S(v)}\varepsilon(H),\qquad
\nu_H(v)+\nu_S(v)=v-1.
\]
As $n-(d+1)-2=q-1$, the composite
\[
\psi_F:\rH^d(D_F;\kk)\xrightarrow{\ \kappa_F\ }H^{d+1}(2^{W_F},D_F;\kk)
\xrightarrow{\ (E_F)_*\ }\Kc_q(F;\Delta)
\]
is an isomorphism. If $F\notin\Delta$ both groups vanish, and $F=N$ contributes nothing since $\Delta\ne2^N$.

We check compatibility with the face-index maps. For $F\subseteq G\subsetneq N$, restriction gives
\[
R_{F,G}:C^\bullet(2^{W_F},D_F)\longrightarrow C^\bullet(2^{W_G},D_G),\qquad
R_{F,G}(e_H^*)=\begin{cases}e_H^*,&H\subseteq W_G,\\0,&\text{otherwise}.\end{cases}
\]

Denote by $H^{d+1}(R_{F,G})$ the induced map
\[
H^{d+1}(R_{F,G}):
H^{d+1}(2^{W_F},D_F;\kk)
\longrightarrow
H^{d+1}(2^{W_G},D_G;\kk).
\]
On costar chains $P_{F,G}^\Delta$ retains $[S]$ exactly when $G\subseteq S$. Since $H\subseteq W_G$ if and only if $G\subseteq N\setminus H$, and $\varepsilon(H)$ does not depend on $F$, we obtain
\[
E_G R_{F,G}=P_{F,G}^\Delta E_F.
\]
The inclusion $D_G\subseteq D_F$ also induces
\[
\operatorname{res}_{F,G}:\rH^d(D_F;\kk)\longrightarrow\rH^d(D_G;\kk).
\]
Naturality of the connecting maps gives $\kappa_G\operatorname{res}_{F,G}=H^{d+1}(R_{F,G})\kappa_F$, and hence
\[
\psi_G\operatorname{res}_{F,G}=\rho^\Delta_{F,G}\psi_F.
\]
The direct sums $\Psi^p=\bigoplus_{|F|=p}\psi_F$ therefore satisfy $\Psi^{p+1}d_{\Lc}=d_{\Kc}\Psi^p$, since the face-index incidence signs agree. This proves $\Lc_q^\bullet(\Delta)\cong\Kc_q^\bullet(\Delta)$.

For $\Delta'\subseteq\Delta''$, we have $D_F(\Delta'')\subseteq D_F(\Delta')$. The corresponding map of relative cochain complexes becomes, under $E_F$, the inclusion of the costar relative chain complexes: both retain the same oriented faces. By naturality of the connecting homomorphisms the isomorphism of complexes commutes with the maps induced by $\Delta'\subseteq\Delta''$. The asserted formula follows from Theorem~\ref{thm:lyub-complex}.

If $\Delta=2^N$, every proper-face link is a nonempty simplex and has zero reduced homology. The only nonzero summand occurs for $F=N$ and $q=n$, where $\rH_{-1}(\{\varnothing\};\kk)=\kk$. Thus $\lambda_{n,n}=1$ and all other Lyubeznik numbers vanish.
\end{proof}

\begin{remark}\label{rem:squarefree-face-link}
There is also an algebraic interpretation. Put $\omega_S=S(-\one_{[n]})$ and $K_q(\Delta)=\Ext_S^{n-q}(\kk[\Delta],\omega_S)$. Graded local duality and the relative form of Hochster's formula give
\[
[K_q(\Delta)]_{\one_F}\cong\operatorname{Hom}_\kk\bigl([\Hloc{q}{\mm}(\kk[\Delta])]_{-\one_F},\kk\bigr)
\cong\Kc_q(F;\Delta).
\]
Multiplication between these squarefree graded pieces corresponds to $\rho^\Delta_{F,G}$, and the associated cellular complex is $\Kc_q^\bullet(\Delta)$. This relates the construction to the double-$\Ext$ formulas and linear strands of Alexander-dual ideals in \cite{Yanagawa2000,AlvarezMontanerVahidi,AlvarezMontanerYanagawa}, without requiring a choice of minimal resolution.
\end{remark}

\subsection{Persistence through costars}\label{ssec:costar-persistence}

The costar complex is functorial also for the full simplex, which was treated separately in Definition~\ref{def:lyub-complex}. For $s\le t$ let $\Phi_q^{s,t}:\Kc_q^\bullet(\Delta^s)\to\Kc_q^\bullet(\Delta^t)$ be the map of Lemma~\ref{lem:face-link-differential}.

\begin{corollary}\label{cor:costar-persistent}
For $s\le t$, with $\Delta^s,\Delta^t\ne2^{[n]}$,
\[
\lambda_{p,q}^{s\to t}=\rank\Bigl(H^p(\Kc_q^\bullet(\Delta^s))
\xrightarrow{\ H^p(\Phi_q^{s,t})\ }H^p(\Kc_q^\bullet(\Delta^t))\Bigr).
\]
The formula also holds when the target is the full simplex. These persistence modules have the finite barcodes and bottleneck stability bound of Theorem~\ref{thm:lyub-stability}.
\end{corollary}

\begin{proof}
For proper levels the formula follows from Theorems~\ref{thm:link-model} and~\ref{thm:lyub-persistent}, including naturality. If $\Delta^t=2^{[n]}$, the only nonzero target group is in bidegree $(n,n)$, and the source is zero in this bidegree unless it too is the full simplex. Thus the map from a proper level to the full simplex is zero, both for Lyubeznik persistence and for the costar complex. When both levels are full simplices the map is the identity.
\end{proof}

For $r\le s\le t$, factoring the structure map through level $s$ gives
\[
0\le\lambda_{p,q}^{r\to t}\le\min\{\lambda_{p,q}(\kk[\Delta^r]),\lambda_{p,q}(\kk[\Delta^t])\},\qquad
\lambda_{p,q}^{r\to t}\le\min\{\lambda_{p,q}^{r\to s},\lambda_{p,q}^{s\to t}\}.
\]
The array $(\lambda_{p,q}^{s\to t})_{0\le p,q\le n}$ is the \emph{persistent Lyubeznik table} from $s$ to $t$; for $s=t$ it is the ordinary Lyubeznik table of $\kk[\Delta^t]$.

\begin{example}\label{ex:additional-lyubeznik-bars}
Let $\Delta$ have facets $\{1,2,3\}$ and $\{1,4,5\}$ on $[5]$. Give the faces of the first triangle filtration value $0$ and all remaining faces value $1$. Then
\[
\Delta^t=\begin{cases}
\text{void},&t<0,\\
2^{\{1,2,3\}},&0\le t<1,\\
2^{\{1,2,3\}}\cup2^{\{1,4,5\}},&t\ge1.
\end{cases}
\]
Both nonvoid levels are contractible, since the two filled triangles meet only at vertex $1$. In particular, $\rH_i(\Delta^t;\kk)=0$ for all $i$ and $t\ge0$.

For $q=3$, every edge link is a point and hence $\Kc_3^2(\Delta^t)=0$. Each triangular facet has link $\{\varnothing\}$ and contributes a copy of $\kk$ to $\Kc_3^3(\Delta^t)$; there are no faces of size $4$. Thus
\[
H^3(\Kc_3^\bullet(\Delta^t))\cong\begin{cases}
0,&t<0,\\
\kk,&0\le t<1,\\
\kk^2,&t\ge1.
\end{cases}
\]
The inclusion from the first nonvoid level to the last sends the class of $\{1,2,3\}$ to the same class, and has rank $1$. Hence $\mathcal B_{3,3}(f)=\{[0,\infty),[1,\infty)\}$.

There is also a bar in bidegree $(1,2)$. Before $t=1$ the links of the vertices present are connected. At $t=1$ the link of vertex $1$ becomes the disjoint union of the edges $\{2,3\}$ and $\{4,5\}$. Hence $\Kc_2^1(\Delta^t)\cong\kk$ for $t\ge1$, and the other vertex summands vanish. Since $\Kc_2^0(\Delta^t)=\rH_1(\Delta^t)=0$ and $\Kc_2^2(\Delta^t)=0$, we get $\lambda_{1,2}=1$ for $t\ge1$ and $0$ before then, so $\mathcal B_{1,2}(f)=\{[1,\infty)\}$.

Neither bar born at $t=1$ occurs in the ordinary persistent reduced homology of this filtration. The example shows that the Lyubeznik modules can record changes in the links of faces even when the homology of the entire complex remains trivial.
\end{example}

\section*{Acknowledgment}

The work of KH was supported in part by the University of Georgia. 
The work of FS was supported in part by the King Fahd University of Petroleum and Minerals. 
The work of GWW was supported in part by the University of Georgia, Georgia Research Alliance, and   NIH R35GM148196.

\section*{AI Usage Disclosure}
The authors used Claude (Anthropic) for grammar and \LaTeX{} checking, proofreading,
literature searching, and generating Python and Macaulay2 code to check the examples.

\bibliographystyle{amsalpha}
\bibliography{reference}

\appendix
\newcommand{\appfig}[1]{\includegraphics[width=\linewidth]{#1}}
\begingroup
\emergencystretch=2em
\section{Computed examples}\label{app:computations}

This appendix computes the invariants of this paper for two small examples, five points in the plane and the two isomers of a small molecule, and compares them with persistent homology and with the persistent graded Betti numbers of \cite{SuwayyidWeiGraphs}. Every Betti, Bass and Lyubeznik number reported below was obtained in two ways: from combinatorial formulas, namely the persistent Hochster formula of \cite{SuwayyidWeiGraphs} for the Betti numbers and the formulas of this paper for the Bass and Lyubeznik numbers, and by computations with free resolutions in Macaulay2. The only exceptions are the Bass numbers at the sets $U$ with more than two elements, which appear in Figure~\ref{fig:app-A}b; they were computed from the formulas alone. The second computation uses no simplicial homology. It rests on the definitions, on a theorem of Musta\c{t}\u{a} \cite{Mustata}, and on algebraic statements of this paper: Lemma~\ref{lem:star-loc}(a) and the base change argument of Proposition~\ref{prop:bass-reduction} for the graded pieces of the Bass numbers, and the straightness of $H^{i}_{I_\Delta}(S)$ together with Proposition~\ref{prop:straight-ext} for the Lyubeznik numbers. The Lyubeznik table of each level is checked in addition against the formula of Yanagawa \cite[Theorem~1.1]{Yanagawa2001} and, at the levels other than the full simplex, against that of \`Alvarez Montaner and Vahidi \cite{AlvarezMontanerVahidi}. The two computations agree in every case; see Section~\ref{app:agreement}.

\subsection{What is computed}\label{app:setting}

Let $X=\{x_1,\dots,x_n\}\subset\mathbb{R}^2$. For $t\ge0$ let $\Delta^t$ be the Vietoris--Rips complex of $X$ at scale $t$: a set $\sigma\subseteq[n]$ is a face if $|x_i-x_j|\le t$ for all $i,j\in\sigma$. For $t<0$ let $\Delta^t$ be void. This is the filtration of the monotone function $f(\sigma)=\max_{i,j\in\sigma}|x_i-x_j|$, with $f(\varnothing)=f(\{i\})=0$. Its critical values $0=t_0<t_1<\dots<t_m$ are $0$ and the distinct distances between the points, and $\Delta^{t_m}=2^{[n]}$. We take $k=\mathbb{Q}$. For the two examples below all results are the same over $\mathbb{F}_2$ and over $\mathbb{F}_3$.

The persistence modules below vanish for $t<0$ and are constant on each interval $[t_i,t_{i+1})$, where $t_{m+1}=\infty$. So every bar has the form $[t_i,t_{j+1})$; for the Betti and Bass modules, whose maps go from $t$ to $s$, such a bar is still the set of values $t$ at which the class is present. Such a module is determined by the ranks $r_{i,j}$ of its structure maps between the levels $t_i$ and $t_j$, $i\le j$, whichever direction the maps have: the multiplicity of the bar $[t_i,t_{j+1})$ is
\[
 r_{i,j}-r_{i-1,j}-r_{i,j+1}+r_{i-1,j+1},
\]
with $r_{i,j}=0$ when an index is out of range. For all pairs of critical values $s\le t$ we computed
\begin{itemize}
\item the ranks of $\widetilde{H}_d(\Delta^s;k)\to\widetilde{H}_d(\Delta^t;k)$, that is, persistent homology;
\item the persistent graded Betti numbers $\beta^{s,t}_{i,\mathbf 1_W}$ of \cite{SuwayyidWeiGraphs}, which are the ranks of the maps $\operatorname{Tor}^S_i(k[\Delta^t],k)_{\mathbf 1_W}\to\operatorname{Tor}^S_i(k[\Delta^s],k)_{\mathbf 1_W}$, and their sums $\beta^{s,t}_{i,j}$ over the sets $W$ with $j$ elements (for $W=[n]$, the number $\beta^{s,t}_{i,\mathbf 1_{[n]}}$ is the rank of $\widetilde{H}_{n-i-1}(\Delta^s;k)\to\widetilde{H}_{n-i-1}(\Delta^t;k)$, by the persistent Hochster formula);
\item the persistent Bass numbers $\mu^{p,s\to t}_{U,G}$ of Definition~\ref{def:pers-bass}, for all $U\subseteq[n]$ and $G\subseteq U^c$ (from the formulas for all $U$; with free resolutions for the sets $U$ with at most two elements, see Section~\ref{app:agreement});
\item the persistent Lyubeznik numbers $\lambda^{s\to t}_{p,q}$ of Definition~\ref{def:pers-lyub}.
\end{itemize}

\subsubsection*{Computation from the formulas}
Persistent homology is computed from the simplicial chain complexes of the $\Delta^{t_i}$. The Betti numbers are computed by the persistent Hochster formula \cite[Theorem~3.2]{SuwayyidWeiGraphs}, the Bass numbers by Theorem~\ref{thm:bass-persistent}, and the Lyubeznik numbers by Corollary~\ref{cor:costar-persistent} with the costar complexes $\mathcal K^\bullet_q(\Delta^t)$ and, for the levels other than the full simplex, a second time by Theorem~\ref{thm:lyub-persistent} with the complexes $\mathcal L^\bullet_q(\Delta^t)$. These computations use only reduced simplicial homology and cohomology of subcomplexes of the $\Delta^{t_i}$, of their links and of their Alexander duals, and the relative homology of the pairs $(\Delta^{t_i},\operatorname{cost}_{\Delta^{t_i}}F)$, with the maps induced by inclusions of complexes and of pairs. This is exact linear algebra on matrices whose rows and columns are indexed by faces.

\subsubsection*{Computation with free resolutions}
Write $I_t=I_{\Delta^t}$.

\emph{Betti numbers.} We lift the surjection $k[\Delta^t]\to k[\Delta^s]$ to a $\mathbb{Z}^n$-graded map between the minimal free resolutions over $S$ and tensor with $k$. Since the resolutions are minimal, this gives the maps on $\operatorname{Tor}$. The rank of the part of multidegree $\mathbf 1_W$ in homological degree $i$ is $\beta^{s,t}_{i,\mathbf 1_W}$.

\emph{Bass numbers.} The module $E^p_t=\operatorname{Ext}^p_S(S/\mathfrak{p}_U,k[\Delta^t])$ is finitely generated, $\mathbb{Z}^n$-graded and annihilated by $\mathfrak{p}_U$. Write $\kappa=\kappa(\mathfrak{p}_U)$. Since $\operatorname{Ext}$ commutes with localization, $\operatorname{Ext}^p_{S_{\mathfrak{p}_U}}(\kappa,k[\Delta^t]_{\mathfrak{p}_U})\cong E^p_t\otimes_S\kappa$ naturally, so $\mu^{p,s\to t}(\mathfrak{p}_U)$ is the rank of $E^p_t\otimes_S\kappa\to E^p_s\otimes_S\kappa$. The matrix of a degree-preserving map between $\mathbb{Z}^n$-graded free $S$-modules has entries $c_{ij}x^{\alpha_j-\beta_i}$, where $\alpha_j$ and $\beta_i$ are the degrees of the basis elements. Over $\kappa$ it becomes $D_1CD_2$, where $D_1$ and $D_2$ are invertible diagonal matrices of monomials in the $x_u$, $u\in U$, and $C$ is the matrix of scalars obtained by setting $x_u=1$ for $u\in U$ and $x_j=0$ for $j\notin U$. So ranks over $\kappa$ are ranks of such scalar matrices. If $\varphi_s$ is a graded presentation matrix of $E^p_s$ and $\psi$ is a graded matrix of the map $E^p_t\to E^p_s$ on generators, then
\[
 \mu^{p,s\to t}(\mathfrak{p}_U)=\operatorname{rank}\bigl[\,C_\psi\mid C_{\varphi_s}\bigr]-\operatorname{rank}C_{\varphi_s}.
\]
These scalar matrices are block diagonal with respect to the part $b\in\mathbb{Z}^{U^c}$ of the degrees in the variables outside $U$. By Lemma~\ref{lem:star-loc}(a) and flat base change, as in the proof of Proposition~\ref{prop:bass-reduction}, the block of $b$ computes the part of degree $b$ of $\operatorname{Ext}^p_{S_{U^c}}(k,k[\Delta^t]_{(U)})$. The block $b=-\mathbf 1_G$ therefore gives $\mu^{p,s\to t}_{U,G}$.

\emph{Lyubeznik numbers.} By \cite[Theorem~1.1]{Mustata} the natural map
\[
 \operatorname{Ext}^i_S(S/I_\Delta,S)\longrightarrow H^{i}_{I_\Delta}(S)
\]
is an isomorphism in all multidegrees $a\ge-\mathbf 1_{[n]}$. It is the map of derived functors induced by the inclusion $\operatorname{Hom}_S(S/I,-)\subseteq\Gamma_I(-)$, so these maps commute with the maps induced on both sides by $I_t\subseteq I_s$, which go from level $s$ to level $t$. The module $H^{i}_{I_\Delta}(S)$ is straight (Section~\ref{sec:lyub-static}), so by Proposition~\ref{prop:straight-ext} the module $\operatorname{Ext}^p_S(k,H^{i}_{I_\Delta}(S))$ is concentrated in multidegree $-\mathbf 1_{[n]}$, and the part of the Koszul complex in this degree involves only graded pieces in multidegrees in $\{-1,0\}^n$. Therefore
\begin{equation}\label{eq:app-lyub-ext-ext}
\begin{aligned}
 \lambda^{s\to t}_{p,q}=\operatorname{rank}\Bigl(&\operatorname{Ext}^p_S\bigl(k,\operatorname{Ext}^{n-q}_S(S/I_s,S)\bigr)_{-\mathbf 1_{[n]}}\\
 &\longrightarrow\operatorname{Ext}^p_S\bigl(k,\operatorname{Ext}^{n-q}_S(S/I_t,S)\bigr)_{-\mathbf 1_{[n]}}\Bigr),
\end{aligned}
\end{equation}
where the map is induced by $S/I_t\to S/I_s$. This holds for all levels, including the full simplex, where $I_t=0$. All modules in \eqref{eq:app-lyub-ext-ext} are finitely generated. So this computation rests on the definition, on Musta\c{t}\u{a}'s theorem, on the straightness of $H^{i}_{I_\Delta}(S)$ and on Proposition~\ref{prop:straight-ext}.

For each level we also evaluated two earlier formulas for the Lyubeznik table. The first is a formula of Yanagawa \cite[Theorem~1.1 and Corollary~3.10]{Yanagawa2001} for the Bass numbers of the modules $H^{j}_{I}(S)$ at the graded maximal ideal $\mathfrak{m}$ of $S$, where $I$ is a squarefree monomial ideal. For $\lambda_{p,q}(S/I)=\mu^p(\mathfrak{m},H^{n-q}_{I}(S))$ it reads
\[
 \lambda_{p,q}(S/I)=\dim_k\bigl[\operatorname{Ext}^{n-p}_S\bigl(\operatorname{Ext}^{n-q}_S(S/I,S),S\bigr)\bigr]_{\mathbf 0}.
\]
The second is the formula of \`Alvarez Montaner and Vahidi \cite[Section~4]{AlvarezMontanerVahidi}, which we use for $I\ne0$: $\lambda_{p,q}(S/I)$ is the dimension of the homology, at the term of homological degree $q-p$, of the $(n-q)$-linear strand of the minimal free resolution of the Alexander dual ideal $I^\vee$ (the resolution of the ideal, with its generators in homological degree $0$), after all variables are set equal to $1$. We use these two formulas only level by level, as a check of the case $s=t$.

\subsection{Example A: five points in the plane}\label{app:exA}

Let $X$ consist of
\[
 x_1=(0,0),\ \ x_2=(-2,1),\ \ x_3=(-2,-1),\ \ x_4=(3,2),\ \ x_5=(3,-2).
\]
These are the vertices of two isosceles triangles of different sizes with the common vertex $x_1$. The critical values are $0$, $2$, $\sqrt5$, $\sqrt{13}$, $4$, $\sqrt{26}$, $\sqrt{34}$. We write a face as the string of its vertices. Figures~\ref{fig:app-A} and~\ref{fig:app-A-betti} show the seven complexes with their facets and, on the scale $t$, the barcodes of the four families of invariants. A bar is an interval $[t_i,t_{j+1})$. Each family has one row, possibly of several lanes, for each index that has at least one bar: the degree $d$, the pair $(p,q)$, the set $U$ (with the colour giving $p$), or the pair $(i,j)$. The number beside a heading counts the bars of the family, and the number beside a row $\beta_{i,j}$ counts the bars of that row. Of the persistent Bass numbers $\mu^{p,s\to t}_{U,G}$ Figure~\ref{fig:app-A}b shows those with $G=\varnothing$ and $U\ne\varnothing$. Those with $U=\varnothing$ are not drawn separately: by Corollary~\ref{cor:bass-betti}(a) they are the persistent graded Betti numbers, $\mu^{p,s\to t}_{\varnothing,G}=\beta^{s,t}_{5-p,\mathbf 1_{[5]\setminus G}}$, so for each $g$ their bars with $|G|=g$ together form the row $\beta_{5-p,5-g}$; for $G=\varnothing$ they are the ranks of $\widetilde H_{p-1}(\Delta^s;k)\to\widetilde H_{p-1}(\Delta^t;k)$, so their bars are the bars of persistent homology in degree $p-1$, which here occur only for $p=1$. So the Bass rows of the figure show only a part of the Bass numbers.

\begin{figure}[htbp]
\centering
\appfig{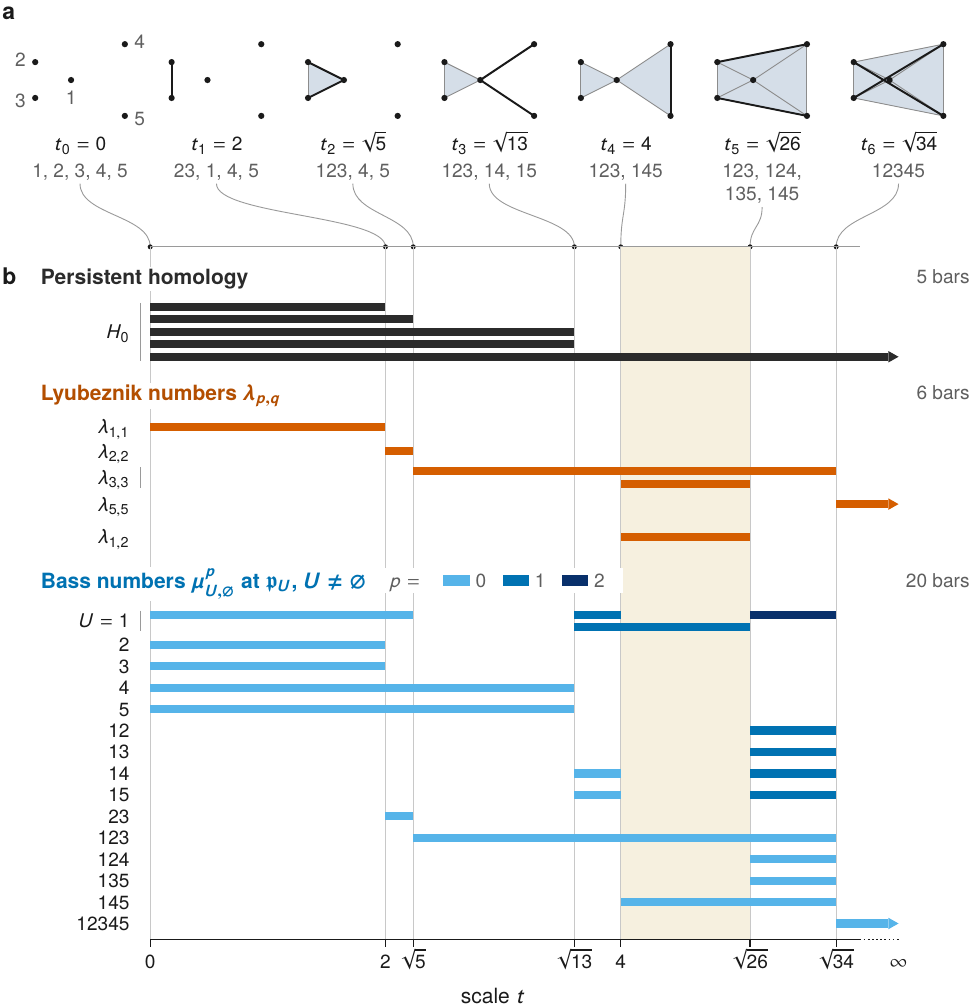}
\caption{Example~A. \textbf{a}, The Rips complexes $\Delta^{t_0}\subset\dots\subset\Delta^{t_6}$ with their facets; the edges that enter at a level are thick. \textbf{b}, The barcodes, on the scale $t$, of persistent homology, of the Lyubeznik modules (Theorem~\ref{thm:lyub-stability}), of the Bass modules $\mathbb M^p_{U,\varnothing}$ at the primes $\mathfrak p_U$ with $U\ne\varnothing$ (Theorem~\ref{thm:bass-stability}; the colour is $p$). An arrow is a bar that does not end, and the numbers at the right count the bars of each family. The band marks $4\le t<\sqrt{26}$, where the bar of $\lambda_{1,2}$ and the second bar of $\lambda_{3,3}$ lie.}\label{fig:app-A}
\end{figure}

\begin{figure}[htbp]
\centering
\appfig{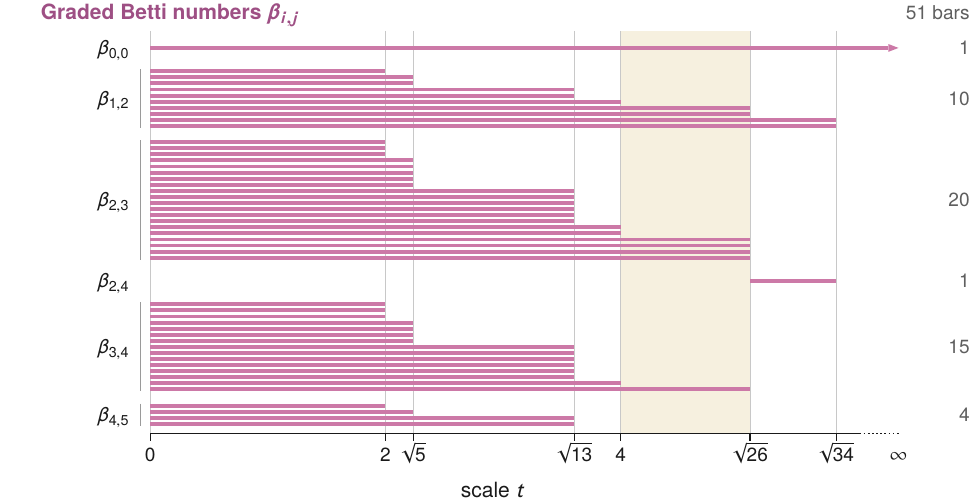}
\caption{Example~A, continued: the barcodes of the persistent graded Betti numbers $\beta_{i,j}$ of \cite{SuwayyidWeiGraphs}, on the scale $t$ and with the band of Figure~\ref{fig:app-A}b. The number beside a row counts its bars.}\label{fig:app-A-betti}
\end{figure}

\begin{figure}[htbp]
\centering
\appfig{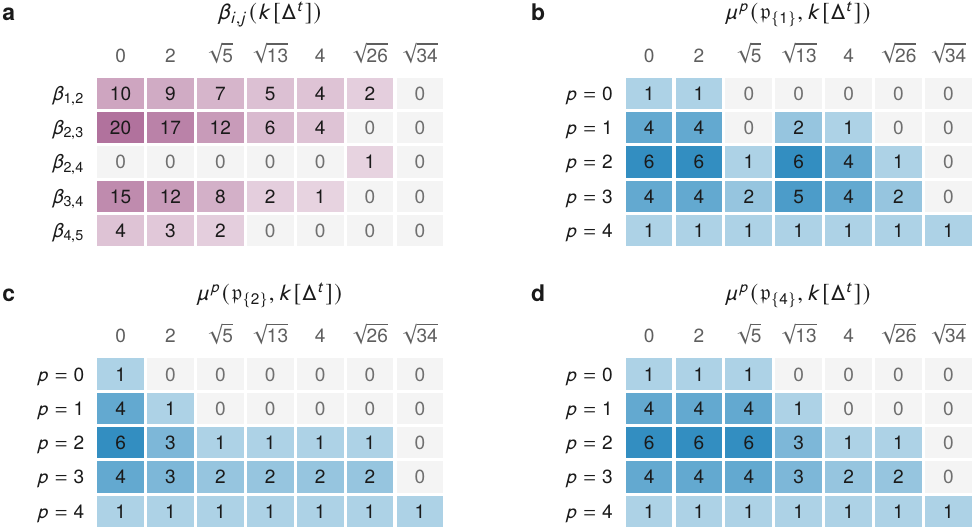}
\caption{Example~A, level by level. \textbf{a}, The graded Betti numbers $\beta_{i,j}(k[\Delta^{t}])$; in addition $\beta_{0,0}=1$ at every level. \textbf{b}--\textbf{d}, The Bass numbers $\mu^p(\mathfrak{p}_{\{v\}},k[\Delta^{t}])$ at the vertex primes for $v=1,2,4$. The columns are the critical values $t$. These are the sums over all $G$, whereas Figure~\ref{fig:app-A}b shows the part with $G=\varnothing$. A darker cell is a larger number.}\label{fig:app-A-numbers}
\end{figure}

\emph{Persistent homology.} The persistent homology has the five bars of Figure~\ref{fig:app-A}b in degree $0$, namely $[0,2)$, $[0,\sqrt5)$ and twice $[0,\sqrt{13})$, and no others. For $t\ge\sqrt{13}$ every facet of $\Delta^t$ contains the vertex $1$, so $\Delta^t$ is a cone and has no reduced homology.

\emph{Graded Betti numbers and Bass numbers at $\mathfrak{m}$.} As noted above, the persistent Bass numbers for $U=\varnothing$ are the persistent graded Betti numbers (Corollary~\ref{cor:bass-betti}(a)). In the computation with free resolutions this says that the maps on $\operatorname{Ext}^p_S(k,-)$ and on $\operatorname{Tor}^S_{5-p}(-,k)$ have the same ranks in corresponding degrees, and this is what we found. The persistent Betti modules have $51$ bars (Figure~\ref{fig:app-A-betti}), and all but one of them start at $0$. These numbers record the combinatorics of the complex: for instance $\beta_{1,2}(k[\Delta^t])$ is the number of pairs of points at distance greater than $t$, which changes at every critical value (Figure~\ref{fig:app-A-numbers}a).

\emph{Bass numbers at vertex primes.} Figures~\ref{fig:app-A-numbers}b--d list $\mu^p(\mathfrak{p}_{\{v\}},k[\Delta^{t_i}])$ for $v=1,2,4$. The reflection $(x,y)\mapsto(x,-y)$ exchanges $2\leftrightarrow3$ and $4\leftrightarrow5$, so the vertices $3$ and $5$ give the same numbers as $2$ and $4$. By Corollary~\ref{cor:bass-betti-link} these are the Betti numbers of $k[\operatorname{lk}_{\Delta^t}\{v\}]$ over the polynomial ring in the other four variables, with $\mu^p=\beta_{4-p}$, and by Corollary~\ref{cor:bass-betti}(b) their persistent versions are the persistent graded Betti numbers of the filtered link of $v$. In degree $p=0$ the barcode at $\mathfrak{p}_{\{v\}}$ is the interval on which $v$ is a facet, that is, an isolated point (Proposition~\ref{prop:bass-bars}):
\[
 [0,\sqrt5)\ \text{for } v=1,\qquad [0,2)\ \text{for } v=2,\qquad [0,\sqrt{13})\ \text{for } v=4 .
\]
In degree $p=1$ with $G=\varnothing$ the module is $\widetilde{H}^0(\operatorname{lk}_{\Delta^t}\{v\})$. It is zero for $v\ne1$, and for $v=1$ its barcode is
\[
 [\sqrt{13},4),\qquad[\sqrt{13},\sqrt{26}) .
\]
The link of the vertex $1$ consists of the edge $23$ and the points $4$ and $5$ at $t=\sqrt{13}$, of the two edges $23$ and $45$ at $t=4$, and it becomes the $4$-cycle with edges $23,35,45,24$ at $t=\sqrt{26}$.

\emph{Lyubeznik numbers.} The Lyubeznik barcodes have $6$ bars in total. On the diagonal $p=q$ they are $[0,2)$ for $q=1$, $[2,\sqrt5)$ for $q=2$, the bars $[\sqrt5,\sqrt{34})$ and $[4,\sqrt{26})$ for $q=3$, and $[\sqrt{34},\infty)$ for $q=5$. As Proposition~\ref{prop:top-row} predicts, every left endpoint below $\sqrt{34}$ is the value of $f$ at a face with $q$ vertices; for $q=3$ these values are $f(123)=\sqrt5$ and $f(145)=4$. The bar $[\sqrt{34},\infty)$ starts at the full simplex, where $\lambda_{5,5}=1$ by Theorem~\ref{thm:lyub-complex}. By Corollary~\ref{cor:lyub-terms}(d), for $t<\sqrt{34}$ the number of bars of $\mathcal B_{q,q}$ that contain $t$ is at most the number of facets of $\Delta^t$ with $q$ vertices, and it can be smaller: $\Delta^{t_5}$ has four facets with three vertices and $\lambda_{3,3}=1$. For $4\le t<\sqrt{26}$ the complex is the union of the triangles $123$ and $145$, which meet only in the vertex $1$. This is the complex of Example~\ref{ex:additional-lyubeznik-bars}, and $\lambda_{1,2}=1$, $\lambda_{3,3}=2$. At $t=\sqrt{26}$ the edges $24$ and $35$ enter, the four triangles $123$, $124$, $145$, $135$ form a disc around the vertex $1$, the link of $1$ becomes connected, and both the bar of $\mathcal B_{1,2}$ and the second bar of $\mathcal B_{3,3}$ end.

\emph{Comparison.} The two bars $[4,\sqrt{26})$ of $\mathcal B_{1,2}$ and $\mathcal B_{3,3}$ are not bars of persistent homology, which vanishes for $t\ge\sqrt{13}$ (i.e. only the infinite bar is left), and no bar of the persistent Betti modules starts at $t=4$. The Bass numbers at $\mathfrak{p}_{\{1\}}$ and the number $\lambda_{1,2}$ both look at the link of the vertex $1$, but they are not the same: for $\sqrt{13}\le t<4$ the link of $1$ has three components and $\mu^1_{\{1\},\varnothing}=2$, while $\lambda_{1,2}=0$. In terms of the costar complex of Definition~\ref{def:face-link-complex}, whose terms are reduced homology groups of links by Lemma~\ref{lem:costar-link}, the links of the other vertices are connected, so $\mathcal K^1_2(\Delta^t)\cong\widetilde{H}_0(\operatorname{lk}_{\Delta^t}\{1\})\cong k^2$, and the differential maps this space isomorphically onto $\mathcal K^2_2(\Delta^t)\cong k^2$, whose two summands come from the facets $14$ and $15$. So here the Bass numbers at $\mathfrak{p}_{\{1\}}$ change whenever the link of $1$ changes, while $\lambda_{1,2}$ keeps only the part of $\widetilde{H}_0(\operatorname{lk}_{\Delta^t}\{1\})$ that is not cancelled by the edges through $1$ that are facets: these are $14$ and $15$ for $\sqrt{13}\le t<4$, and there are none for $4\le t<\sqrt{26}$.

Since $\beta_{1,\mathbf 1_W}(k[\Delta])\ne0$ exactly when $W$ is a minimal nonface, the multigraded Betti numbers determine $\Delta^t$, and so in principle they determine every other invariant considered here. The comparison concerns what each invariant displays. The Betti barcodes of Example~A have $51$ bars, and the Betti numbers change at every critical value. The Lyubeznik barcodes have $6$ bars, and by \cite{AlvarezMontanerYanagawa} the Lyubeznik table of each level depends only on the homeomorphism type of $|\Delta^t|$ and on the characteristic of $k$.

\subsection{Example B: two isomers}\label{app:exB}

Section~6.2 of \cite{SuwayyidWeiGraphs} applies persistent Betti numbers to two isomers, the \emph{cis} and the \emph{trans} form of 1,2-dichloroethene ($\mathrm{C_2H_2Cl_2}$). We take the same two molecules, each as a set of six points, in an idealized planar geometry. The bond lengths are $1.33$ for the double bond between the carbon atoms, $1.73$ for a carbon--chlorine bond and $1.09$ for a carbon--hydrogen bond, in \AA, and all bond angles are $120^\circ$. Explicitly,
\begin{gather*}
 \mathrm C_1=(-0.665,0),\qquad \mathrm C_2=(0.665,0),\\
 \mathrm{Cl}_1=\mathrm C_1+1.73\,(\cos120^\circ,\sin120^\circ),\\
 \mathrm H_1=\mathrm C_1+1.09\,(\cos240^\circ,\sin240^\circ),
\end{gather*}
and $\mathrm{Cl}_2$, $\mathrm H_2$ are the images of $\mathrm{Cl}_1$, $\mathrm H_1$ under the reflection $(x,y)\mapsto(-x,y)$ for the \emph{cis} isomer and under the point reflection $(x,y)\mapsto(-x,-y)$ for the \emph{trans} isomer. The atoms $\mathrm C_1,\mathrm C_2,\mathrm{Cl}_1,\mathrm{Cl}_2,\mathrm H_1,\mathrm H_2$ are the vertices $1,\dots,6$. All squared distances are rational, so the two Rips filtrations are computed exactly; we print critical values to three decimals.

The two filtrations coincide for $t<2.420$: the edges $15$ and $26$ enter at $1.090$, $12$ at $1.330$, $13$ and $24$ at $1.730$, and $16$ and $25$ at $2.099$.
The remaining pairs enter as follows: for the \emph{cis} isomer $56$ at $2.420$, $35$ and $46$ at $2.463$, $14$ and $23$ at $2.658$, $34$ at $3.060$, $36$ and $45$ at $3.670$; for the \emph{trans} isomer $35$ and $46$ at $2.463$, $14$ and $23$ at $2.658$, $36$ and $45$ at $2.795$, $56$ at $3.069$, $34$ at $4.283$.
Figure~\ref{fig:app-B-molecules} shows the two molecules and the two complexes for $2.658\le t<2.795$, and Figure~\ref{fig:app-B} shows the two filtrations and the barcodes of persistent homology, of the Lyubeznik modules and of the Bass modules at the double bond. Figure~\ref{fig:app-B-betti} shows the barcodes of the persistent graded Betti numbers. In both figures carbon is dark, chlorine green and hydrogen white, and in Figure~\ref{fig:app-B} the edges that enter at a level are thick. In each row the bars of the \emph{cis} isomer come first, and an empty lane means that the isomer has no bar in that row. A vertical line through all rows is a critical value of both filtrations; a line that runs only through the lanes of one isomer is a critical value of that isomer alone. The values $3.060$ (\emph{cis}) and $3.069$ (\emph{trans}) are so close that their two lines together look like a single broken line through the lanes of both isomers; each is a critical value of one isomer only.

\begin{figure}[htbp]
\centering
\appfig{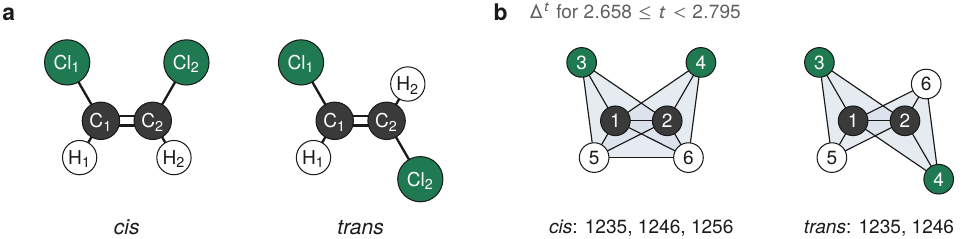}
\caption{Example~B. \textbf{a}, The two isomers in the idealized geometry. \textbf{b}, The complexes $\Delta^t$ for $2.658\le t<2.795$ with their facets; the vertices $1,\dots,6$ are $\mathrm C_1,\mathrm C_2,\mathrm{Cl}_1,\mathrm{Cl}_2,\mathrm H_1,\mathrm H_2$.}\label{fig:app-B-molecules}
\end{figure}

\begin{figure}[htbp]
\centering
\appfig{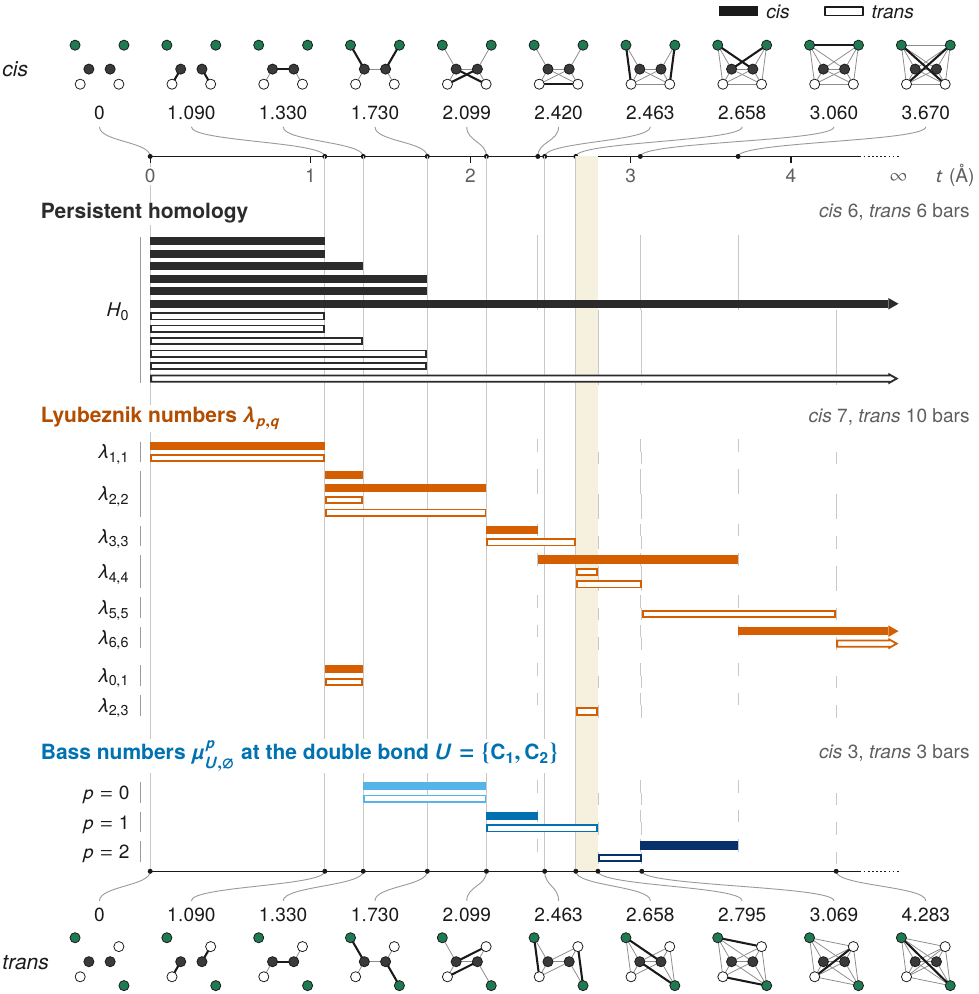}
\caption{Example~B: the two Rips filtrations and the barcodes. Above (\emph{cis}) and below (\emph{trans}): the graph of each level, with its new edges thick, joined by a curve to its critical value on the scale $t$. Between them: the barcodes of persistent homology, of the Lyubeznik modules, and of the Bass modules $\mathbb M^p_{U,\varnothing}$ at the double bond $U=\{\mathrm C_1,\mathrm C_2\}$, with the \emph{cis} lanes (filled bars) above the \emph{trans} lanes (open bars). An arrow is a bar that does not end; the numbers at the right count bars. A broken vertical line is a critical value of one isomer alone; those at $3.060$ and $3.069$ nearly coincide. The band marks $2.658\le t<2.795$, the scales of Figure~\ref{fig:app-B-molecules}b.}\label{fig:app-B}
\end{figure}

\begin{figure}[htbp]
\centering
\appfig{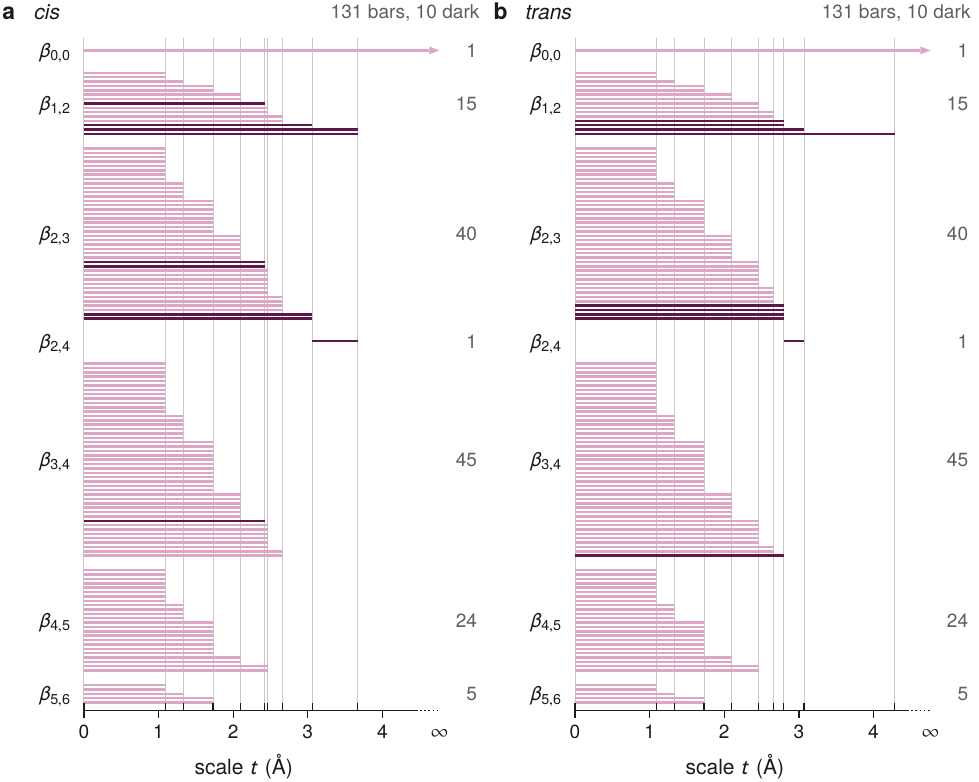}
\caption{Example~B: the barcodes of the persistent graded Betti numbers $\beta_{i,j}$ of \cite{SuwayyidWeiGraphs}; \textbf{a}, the \emph{cis} isomer; \textbf{b}, the \emph{trans} isomer. A bar is dark if the other isomer has no bar with the same ends in the same row, bars being counted with multiplicity; the two barcodes share $121$ of their $131$ bars. The vertical lines mark the critical values of each filtration, which are those of Figure~\ref{fig:app-B}, the number beside a row counts its bars, and an arrow is a bar that does not end.}\label{fig:app-B-betti}
\end{figure}

\emph{Persistent homology} does not separate the two isomers. In reduced homology of degree $0$ both have the bars $[0,1.090)$ twice, $[0,1.330)$ once and $[0,1.730)$ twice, and neither has a bar in a higher degree.

\emph{The Lyubeznik barcodes differ.} Both isomers have the bar $[1.090,1.330)$ of $\mathcal B_{0,1}$. Apart from it, the only bar off the diagonal $p=q$ is the bar $[2.658,2.795)$ of $\mathcal B_{2,3}$ of the \emph{trans} isomer. For $2.658\le t<2.795$ the \emph{trans} complex is the union of the two tetrahedra $1235$ and $1246$, that is, $\{\mathrm C_1,\mathrm C_2,\mathrm{Cl}_1,\mathrm H_1\}$ and $\{\mathrm C_1,\mathrm C_2,\mathrm{Cl}_2,\mathrm H_2\}$, which meet only in the edge $\mathrm C_1\mathrm C_2$ (Figure~\ref{fig:app-B-molecules}b). The link of this edge is the disjoint union of the edges $35$ and $46$. For $q=3$ the only nonzero term of the costar complex is $\mathcal K^2_3(\Delta^t)\cong\widetilde{H}_0(\operatorname{lk}_{\Delta^t}\{1,2\})\cong k$, which gives $\lambda_{2,3}=1$, and for $q=4$ the only nonzero term is $\mathcal K^4_4(\Delta^t)\cong k^2$, from the two facets, which gives $\lambda_{4,4}=2$. This complex is the cone with apex $1$ over the two triangles $235$ and $246$, which meet in a vertex as in Example~\ref{ex:additional-lyubeznik-bars}. For the \emph{cis} isomer at the same scales the atoms $\mathrm H_1$ and $\mathrm H_2$ are already joined, the complex consists of the three tetrahedra $1235$, $1256$, $1246$, and consecutive ones meet in a triangle. This is a shellable $3$-ball, so $k[\Delta^t]$ is Cohen--Macaulay \cite{BrunsHerzog} and its Lyubeznik table is trivial \cite{AlvarezMontanerYanagawa}: the only nonzero Lyubeznik number is $\lambda_{4,4}=1$.

\emph{Bass numbers at the double bond.} The same difference is seen at the prime $\mathfrak{p}_U$ of $U=\{\mathrm C_1,\mathrm C_2\}$. The module $\mathbb M^1_{U,\varnothing}$ is $\widetilde{H}^0$ of the link of the edge $\mathrm C_1\mathrm C_2$ (Theorem~\ref{thm:bass-persistent}), and its bar is $[2.099,2.420)$ for the \emph{cis} isomer and $[2.099,2.795)$ for the \emph{trans} isomer (Figure~\ref{fig:app-B}): the link of the double bond stays disconnected for longer in the \emph{trans} isomer.

\emph{Graded Betti numbers.} The persistent graded Betti numbers also separate the two filtrations (Figure~\ref{fig:app-B-betti}). For instance $\beta_{1,2}(k[\Delta^t])$ is the number of pairs of atoms at distance greater than $t$. The barcodes of the modules $\beta_{i,j}$ have $131$ bars for the \emph{cis} isomer and $131$ for the \emph{trans} isomer, and $121$ bars of each, counted with multiplicity, also occur in the other, in the same row and with the same ends; the Lyubeznik barcodes have $7$ and $10$ bars. Section~6.2 of \cite{SuwayyidWeiGraphs} compares the persistent Betti numbers of the two isomers qualitatively and gives no numerical values. The setting here, with the rings $k[\Delta^t]$ of the Rips complexes and an idealized geometry, need not coincide with the one used there, so the two are not compared quantitatively.

\subsection{Agreement of the two computations}\label{app:agreement}

Figure~\ref{fig:app-agreement} lists the filtrations for which both computations of Section~\ref{app:setting} were carried out: Examples~A and~B; the filtration of Example~\ref{ex:additional-lyubeznik-bars}, a triangle to which a second triangle is attached at a vertex; the six-vertex triangulation of the real projective plane of Example~\ref{ex:rp2}, entering in two steps, first its vertices and edges and then its triangles; the Rips filtrations of further planar point sets, namely the vertices of a regular hexagon, a square with a triangle attached at a vertex, two unit squares with a common vertex, and two sets of random points with integer coordinates; and three filtrations generated by random faces, some of whose levels are not flag complexes. Figure~\ref{fig:app-agreement}a gives the numbers of nonzero ranks $\beta^{s,t}_{i,\mathbf 1_W}$, $\mu^{p,s\to t}_{U,G}$ and $\lambda^{s\to t}_{p,q}$, over all pairs $s\le t$ of critical values, that were computed in both ways. For the Bass numbers only the sets $U$ with at most two elements are counted (at most one element for the entry marked $\ast$).
The point sets are the following. Regular hexagon: $(\pm2,0)$, $(\pm1,\pm\sqrt3)$. Square and triangle: $(0,0)$, $(4,0)$, $(4,4)$, $(0,4)$, $(-3,-1)$, $(-1,-3)$. Two squares: $(0,0)$, $(1,0)$, $(1,1)$, $(0,1)$, $(-1,0)$, $(-1,-1)$, $(0,-1)$. Random points I: $(4,18)$, $(14,15)$, $(3,15)$, $(20,12)$, $(2,8)$, $(6,3)$. Random points II: $(18,2)$, $(19,0)$, $(7,18)$, $(11,19)$, $(15,20)$, $(17,4)$. The last three filtrations are not Rips filtrations: all six vertices are present at time $0$, and at time $j=1,\dots,7$ the $j$-th face of the following list enters together with its subsets. Random faces I: $1236,\ 3456,\ 14,\ 1346,\ 2456,\ 145,\ 1256$. Random faces II: $1235,\ 46,\ 156,\ 1234,\ 1456,\ 1256,\ 2356$. Random faces III: $26,\ 136,\ 256,\ 34,\ 346,\ 345,\ 24$.

\begin{figure}[htbp]
\centering
\appfig{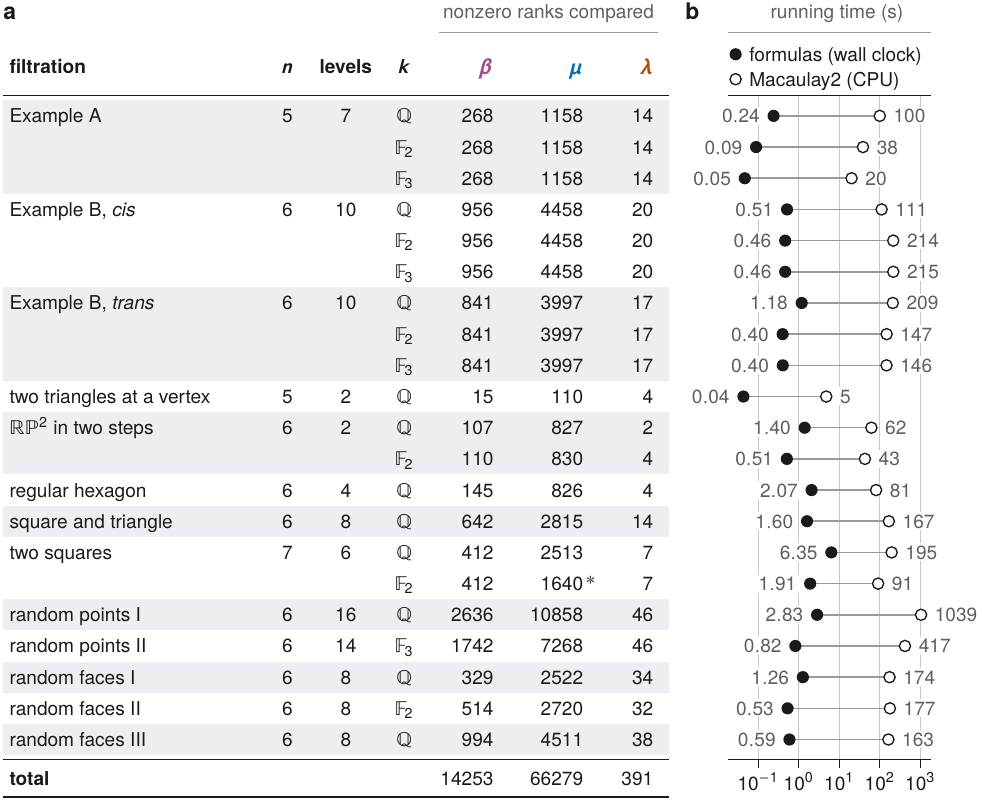}
\caption{Agreement of the two computations. \textbf{a}, For each filtration and field $k$ ($n$ is the number of vertices, and levels is the number of critical values): the numbers of nonzero ranks $\beta^{s,t}_{i,\mathbf 1_W}$, $\mu^{p,s\to t}_{U,G}$ (for $|U|\le2$; for $|U|\le1$ in the entry marked $\ast$) and $\lambda^{s\to t}_{p,q}$, over all pairs $s\le t$ of critical values, that were computed both from the formulas and with free resolutions in Macaulay2. The two computations give the same value in every case. \textbf{b}, Running times in seconds on a logarithmic scale: wall-clock time for the formulas and CPU time for Macaulay2 (see the text).}\label{fig:app-agreement}
\end{figure}

All $80923$ pairs of numbers agree. At every level of every filtration the Lyubeznik table also agrees with Yanagawa's formula, and at every level other than the full simplex it agrees with the formula of \`Alvarez Montaner and Vahidi. For the projective plane the two computations agree over $\mathbb{Q}$ and over $\mathbb{F}_2$, where the Lyubeznik tables differ as in Example~\ref{ex:rp2}.

Figure~\ref{fig:app-agreement}b gives running times in seconds. For the formulas this is the wall-clock time of a direct Python implementation with exact arithmetic; for the computation with free resolutions it is the CPU time of Macaulay2 without start-up. For each row both computations ran, one after the other, on one core of the same machine. Different rows were run as separate jobs, on machines that need not be the same, so times should be compared only within a row. The first computation includes the Bass numbers at all $2^n$ sets $U$ and both complexes $\mathcal K^\bullet_q$ and $\mathcal L^\bullet_q$. The second computes Bass numbers only at the sets $U$ with at most two elements (at most one for the entry marked $\ast$), includes the formulas of Yanagawa and of \`Alvarez Montaner and Vahidi, and is dominated by the modules $\operatorname{Ext}^p_S(S/\mathfrak{p}_U,k[\Delta^t])$ and the maps between them.

\endgroup

\end{document}